\documentclass{amsart}
\usepackage{amssymb}
\usepackage[all]{xy}
\CompileMatrices
\newtheorem*{introthm}{Theorem}
\newtheorem*{introprop}{Proposition}
\newtheorem{thm}{Theorem}[section]
\newtheorem{lemma}[thm]{Lemma}
\newtheorem{prop}[thm]{Proposition}
\newtheorem{coro}[thm]{Corollary}
\theoremstyle{definition}
\newtheorem{defi}[thm]{Definition}
\newtheorem{exam}[thm]{Example}
\newtheorem{rem}[thm]{Remark}

\numberwithin{equation}{thm}
\def\nzd{{\rm nzd}}

\def\mal{\! \cdot \!}

\def\rq#1{\widehat{#1}}
\def\t#1{\widetilde{#1}}
\def\b#1{\overline{#1}}
\def\bangle#1{\langle #1 \rangle}

\def\KK{{\mathbb K}}

\def\ZZ{{\mathbb Z}}

\def\QQ{{\mathbb Q}}
\def\PP{{\mathbb P}}
\def\Of{{\mathcal{O}}}

\def\Pic{\operatorname{Pic}}
\def\Piclin{\Pic_{\rm lin}}
\def\Chara{\operatorname{Char}}
\def\Hom{{\rm Hom}}
\def\Mor{{\rm Mor}}

\def\Spec{{\rm Spec}}

\begin{document}
\title[Homogeneous coordinates]
      {Homogeneous coordinates \\ 
       for algebraic varieties}
\author[F. Berchtold]{Florian Berchtold} 
\address{Fachbereich Mathematik und Statistik, Universit\"at Konstanz}
\email{berchtof@fmi.uni-konstanz.de}
\author[J.~Hausen]{J\"urgen Hausen} 
\address{Fachbereich Mathematik und Statistik, Universit\"at Konstanz}
\email{Juergen.Hausen@uni-konstanz.de}
\begin{abstract}
We associate to every divisorial (e.g.~smooth) variety $X$ with
only constant invertible global functions and finitely generated 
Picard group a $\Pic(X)$-graded homogeneous coordinate ring. 
This generalizes the usual homogeneous coordinate ring of the
projective space and constructions of Cox and Kajiwara for smooth 
and divisorial toric varieties.
We show that the homogeneous coordinate ring defines in fact 
a fully faithful functor.
For normal complex varieties $X$ with only
constant global functions, we even obtain an equivalence of
categories. 
Finally, the homogeneous coordinate ring of a locally factorial 
complete irreducible variety with free finitely generated
Picard group turns out to be a Krull ring admitting 
unique factorization.
\end{abstract}

\maketitle

\section*{Introduction}

The principal use of homogeneous coordinates is 
that they relate
the geometry of algebraic varieties to the theory 
of graded rings. 
The classical example is the projective $n$-space:
its homogeneous coordinate ring is the polynomial 
ring in $n+1$ variables, graded by the usual degree.
Cox~\cite{Co} and Kajiwara~\cite{Ka} 
introduced homogeneous coordinate rings for 
toric varieties.
Cox's construction is meanwhile a
standard instrument in toric geometry;
for example, it is used in~\cite{BrVe} to prove
an equivariant Riemann-Roch Theorem, and in~\cite{MuSmTs} 
for a description of $\mathcal{D}$-modules on toric 
varieties. 

In this article, we construct homogeneous 
coordinates for a fairly general class of algebraic 
varieties:
Let $X$ be a divisorial variety --- e.g.~$X$ is
$\QQ$-factorial or quasiprojective~\cite{Bo} --- 
such that $X$ has only constant globally invertible 
functions and the Picard group $\Pic(X)$ is finitely 
generated.
If the (algebraically closed) ground field $\KK$ is 
of characteristic $p > 0$, then we require that
the multiplicative group $\KK^{*}$ is of infinite 
rank over $\ZZ$, and that $\Pic(X)$ has no 
$p$-torsion.
Examples of such varieties are complete smooth rational complex 
varieties. Moreover, all Calabi-Yau varieties fit into this
framework. 

To define the homogeneous coordinate ring of $X$, consider 
a family of line bundles $L$ on $X$ such that the classes
$[L]$ generate $\Pic(X)$. 
Choosing a common trivializing cover $\mathfrak{U}$
for the bundles $L$, one can achieve that they form a 
finitely generated free abelian group $\Lambda$, 
which is isomorphic to a subgroup of the group of 
cocycles $H^{1}(\mathcal{O}^{*}, \mathfrak{U})$. 
The sheaves of sections $\mathcal{R}_{L}$, where 
$L \in \Lambda$, then fit together to a sheaf $\mathcal{R}$ 
of $\Lambda$-graded $\mathcal{O}_{X}$-algebras. 
Such sheaves $\mathcal{R}$ and their global sections 
$\mathcal{R}(X)$ are often studied. 
For example, in~\cite{HK} they have been used to 
characterize when Mori's program can be carried out, 
and in~\cite{Ha1} they are the starting point
for quotient constructions in the spirit
of Mumford's Geometric Invariant Theory.

A first important observation is that 
we can pass from the above $\Lambda$-graded 
$\mathcal{O}_{X}$-algebras $\mathcal{R}$ to a 
universal $\mathcal{O}_{X}$-algebra $\mathcal{A}$,
which is graded by the Picard group $\Pic(X)$.
This solves in particular the ambiguity problem 
mentioned in~\cite[Remark p.~341]{HK}.
More precisely, we introduce in Section~\ref{section3} 
the concept of a {\em shifting family\/} for the 
$\mathcal{O}_{X}$-algebra $\mathcal{R}$. 
This enables us to identify in a systematic manner
two homogeneous parts $\mathcal{R}_{L}$ and
$\mathcal{R}_{L'}$ if $L$ and $L'$ define the same 
class in $\Pic(X)$.
The result is a projection $\mathcal{R} \to \mathcal{A}$
onto a $\Pic(X)$-graded $\mathcal{O}_{X}$-algebra $\mathcal{A}$.

The {\em homogeneous coordinate ring\/} of $X$ then 
is a pair $(A,\mathfrak{A})$.
The first part $A$ is the $\Pic(X)$-graded $\KK$-algebra of
global sections $\mathcal{A}(X)$. 
The meaning of the second part $\mathfrak{A}$ is 
roughly speaking the following: 
It turns out that $A$ is the algebra of functions of
a quasiaffine variety $\rq{X}$.
Such algebras need not of finite type over $\KK$, 
and $\mathfrak{A}$ is a datum describing 
all the possible affine closures of $\rq{X}$. 
From the algebraic point of view,
the homogeneous coordinate ring is
a {\em freely graded quasiaffine algebra}; 
the category of such algebras is introduced and 
discussed in Sections~\ref{section1} and~\ref{section2}. 

The first main result of this article is that the 
homogeneous coordinate ring is indeed functorial, 
that means that given a morphism $X \to Y$ of varieties, 
we obtain a morphism of the associated freely graded 
quasiaffine algebras, see Section~\ref{section5}. 
In fact, we prove much more, see
Theorem~\ref{fullyfaithful}: 

\begin{introthm}
The assignment $X \mapsto (A,\mathfrak{A})$ 
is a fully faithful functor 
from the category of divisorial varieties $X$ 
with finitely generated Picard group and 
$\mathcal{O}^{*}(X) = \KK^{*}$ to the category of freely graded
quasiaffine algebras.
\end{introthm}

Note that this statement generalizes in particular the description of
the set $\Hom(X,Y)$ of morphisms of two divisorial toric varieties
$X$, $Y$ obtained by Kajiwara in~\cite[Cor.~4.9]{Ka}. 
In the toric situation, $\mathcal{O}^{*}(X) = \KK^{*}$ is a 
usual nondegeneracy assumption: it just means that $X$ has no torus 
factors.

Having proved Theorem~\ref{fullyfaithful}, the task is to
translate geometric properties of a given variety $X$ 
to algebraic properties of its homogeneous coordinate ring
$(A,\mathfrak{A})$. 
In Section~\ref{section6}, we do this for basic properties
of $X$, like smoothness and normality. 
In the latter case, the $\KK$-algebra $A$ is a normal Krull ring. 
Moreover, we discuss quasicoherent sheaves,
and we give descriptions of affine morphisms and closed
embeddings.

In our second main result, we restrict to normal divisorial 
varieties $X$ with finitely generated Picard group and 
$\mathcal{O}(X) = \KK$.
We call such varieties {\em tame}. 
The homogeneous coordinate ring 
$(A,\mathfrak{A})$ of a tame variety $X$ is 
{\em pointed\/} in the sense that 
$A$ is normal with $A_{0} = \KK$ and $A^{*} = \KK^{*}$.
Moreover, $(A,\mathfrak{A})$ is {\em simple\/} in the sense
that the corresponding quasiaffine variety $\rq{X}$
admits only trivial ``linearizable'' bundles, 
see Section~\ref{section7} for the 
precise definition. 
In Theorem~\ref{equivthm}, we show:  

\begin{introthm}
The assignment $X \mapsto (A,\mathfrak{A})$ 
defines an equivalence of the category of tame varieties with the
category of simple pointed algebras.
\end{introthm}

Specializing further to the case of a free Picard group gives the
class of {\em very tame\/} varieties, see Section~\ref{section8}. 
Examples are the Grassmannians and all smooth complete toric varieties. 
For this class, we obtain a nice description of products in 
terms of homogeneous coordinate rings, see Proposition~\ref{products}.
The possibly most remarkable observation is that very tame varieties
open a geometric approach to unique factorization conditions
for multigraded Krull rings, see Proposition~\ref{freefactorial}: 

\begin{introprop}
A very tame variety is locally factorial if and only if its
homogeneous coordinate ring is a unique factorization domain.
\end{introprop}

We conclude the article with an example 
underlining this principle: 
Let $X$ be the projective line with the points $0,1$ 
and $\infty$ doubled, that means that $X$ is nonseparated. 
Nevertheless, $X$ is very tame and its Picard group is
isomorphic to $\ZZ^{4}$. 
As mentioned before, $A = \mathcal{A}(X)$ 
is a unique factorization domain. 
It turns even out to be a classical example of 
a factorial singularity, namely
$$ 
A 
= 
\KK[T_{1}, \ldots, T_{6}]/\bangle{T_{1}^{2}+ \ldots + T_{6}^{2}}.
$$

The quasiaffine variety $\rq{X}$ corresponding to the homogeneous 
coordinate ring of $X$ is an open subset of $\Spec(A)$. 
The prevariety $X$ is a geometric quotient of $\rq{X}$ by a free action 
of a fourdimensional algebraic torus. In particular, $\rq{X}$
is locally isomorphic to the toric variety $\KK \times (\KK^{*})^{4}$.
That means that $\rq{X}$ is toroidal, even with respect to the 
Zariski Topology, but not toric.
\tableofcontents

\section{Quasiaffine algebras and quasiaffine varieties}\label{section1}

Throughout the whole article we work in the category of algebraic
varieties following the setup of~\cite{Ke}. In particular, we work
over an algebraically closed field $\KK$, and the word point always
refers to a closed point. Note that in our setting a variety is
reduced but it need neither be separated nor irreducible.

The purpose of this section is to provide an algebraic
description of the category of quasiaffine varieties. 
The idea is very simple: Every quasiaffine variety $X$ is an open
subset of an affine variety $X'$ and hence is described by the
inclusion $\mathcal{O}(X') \subset \mathcal{O}(X)$ and the vanishing
ideal of the complement $X' \setminus X$ in
$\mathcal{O}(X')$. 

However, in general the algebra of functions $\mathcal{O}(X)$ of a
quasiaffine variety $X$ is not of finite type, see for example~\cite{Re}.
Thus there is no canonical choice of an affine closure $X'$ for a
given $X$. To overcome this ambiguity, we have to treat all possible
affine closures at once. 

We introduce the necessary algebraic notions. 
By a $\KK$-algebra we always mean a reduced commutative 
algebra $A$ over $\KK$ having a unit element. 
We write $\bangle{I}$ for the ideal generated by a 
subset $I \subset A$.
The set of nonzerodivisors of a $\KK$-algebra $A$ is
denoted by $\nzd(A)$. Recall that we have a canonical inclusion 
$A \subset \nzd(A)^{-1}A$ into the algebra of fractions.

\goodbreak

\begin{defi}\label{closedsubalgebra}
Let $A$ be a $\KK$-algebra.
\begin{enumerate}
\item A {\em closing subalgebra\/} of $A$ is a pair $(A',I')$ where
  $A' \subset A$ is a subalgebra of finite type over $\KK$ and $I'
  \subset A'$ is an ideal in $A'$ with
$$ 
I' = \sqrt{\bangle{I' \cap \nzd(A)}},
\qquad
A = \bigcap_{f \in I' \cap \nzd(A)} A_{f}, 
\qquad 
A'_{f} = A_{f} \text{ for all } f \in I'.
$$ 
\item Two closing subalgebras $(A',I')$ and $(A'',I'')$ of $A$ are
  called {\em equivalent\/} if there is a closing subalgebra
  $(A''',I''')$ of $A$ such that
$$ A' \cup A'' \subset A''', 
   \qquad 
   I''' = \sqrt{\bangle{I'}} = \sqrt{\bangle{I''}}. $$
\end{enumerate}
\end{defi}

Note that \ref{closedsubalgebra}~(ii) does indeed define an
equivalence relation. In terms of these notions, the algebraic data to
describe quasiaffine varieties are the following: 

\begin{defi}\label{quasiaffalgdef}
\begin{enumerate}
\item A {\em quasiaffine algebra\/} is a pair $(A,\mathfrak{A})$,
  where $A$ is a $\KK$-algebra and $\mathfrak{A}$ is the
  equivalence class of a closing subalgebra $(A',I')$ of $A$.
\item A {\em homomorphism\/} of quasiaffine algebras
  $(B,\mathfrak{B})$ and $(A,\mathfrak{A})$ is a homomorphism $\mu
  \colon B \to A$ such that there exist $(B',J') \in \mathfrak{B}$
  and $(A',I') \in \mathfrak{A}$ with
$$ \mu(B') \subset A', \qquad I' \subset \sqrt{\bangle{\mu(J')}}. $$
\end{enumerate}
\end{defi}

We show now that the category of quasiaffine varieties
is equivalent to the category of quasiaffine algebras by associating to
every variety $X$ an equivalence class $\mathfrak{O}(X)$ 
of closing subalgebras of $\mathcal{O}(X)$. We use the following notation: 
Given a variety $X$ and a regular function $f \in \mathcal{O}(X)$, let   
$$ X_{f} := \{x \in X; \; f(x) \ne 0\}.$$

\begin{defi}\label{naturalpairdef}
Let $X$ be a quasiaffine variety. Let 
$A' \subset \mathcal{O}(X)$ be a subalgebra of finite type
and $I' \subset A'$ a radical ideal. 
We call $(A',I')$ a {\em natural pair\/} on $X$, 
if for every $f \in I'$ the set $X_{f}$ is affine with
$\mathcal{O}(X_{f}) = A'_{f}$ and the sets $X_{f}$, $f \in I'$, 
cover $X$. We define $\mathfrak{O}(X)$ to be the collection of all
natural pairs on $X$. 
\end{defi}

So, our first task is to verify that the collection
$\mathfrak{O}(X)$ is in fact an equivalence class of 
closing subalgebras of $\mathcal{O}(X)$. This is done in two 
steps:
 
\begin{lemma}\label{naturalpairs}
Let $X$ be a quasiaffine variety. Let $(A',I')$ be a 
natural pair on $X$, and set $X' := \Spec(A')$. 
\begin{enumerate}
\item The morphism $X \to X'$ defined by $A' \subset \mathcal{O}(X)$ 
  is an open embedding, $I'$ is the vanishing ideal of 
  $X' \setminus X$, and $(A',I')$ is a closing subalgebra of
  $\mathcal{O}(X)$.
\item For a subalgebra $A'' \subset \mathcal{O}(X)$ of finite type
  with $A' \subset A''$, consider the ideal $I'' :=
  \sqrt{\bangle{I'}}$ of $A''$. 
  Then $(A'',I'')$ is a natural pair on $X$. 
\end{enumerate} 
\end{lemma}

\proof Recall that for any $f \in \mathcal{O}(X)$
we have $\mathcal{O}(X_{f}) = \mathcal{O}(X)_{f}$. 
In particular, $X \to X'$ is locally given by 
isomorphisms $X_{f} \to X'_{f}$, $f \in I'$. 
This implies that $X \to X'$ is an open embedding and that
$I' \subset A'$ is the vanishing ideal of $X' \setminus X$.
Finally, $(A',I')$ is a closing subalgebra, because  
up to passing to the radical, $I'$ is generated
by the $f \in I'$ that are nontrivial on each irreducible 
component of $X$.

We turn to assertion~(ii). Let $X'' := \Spec(A'')$. 
It suffices to verify that the morphism $X \to X''$ defined 
by $A'' \subset \mathcal{O}(X)$ is an open embedding and 
that $I'' \subset A''$ is the vanishing ideal of 
the complement $X'' \setminus X$. 
Again this holds, because for every $f \in I'$ the map $X \to X''$ 
restricts to an isomorphism $X_{f} \to X''_{f}$. \endproof

\begin{lemma}\label{quasiaff2closingsubalg}
The collection $\mathfrak{O}(X)$ of all natural pairs 
on a quasiaffine variety $X$ is an equivalence
class of closing subalgebras of $\mathcal{O}(X)$. 
\end{lemma}

\proof First note that there exist natural pairs $(A',I')$ on $X$, 
because for every affine closure $X \subset X'$ we obtain such a pair
by setting $A' := \mathcal{O}(X')$ and defining $I' \subset A'$ to be 
the vanishing ideal of the complement $X' \setminus X$. Moreover, by
Lemma~\ref{naturalpairs}~(i), we know that every natural pair is a closing
subalgebra of $\mathcal{O}(X)$.

We show that any two natural pairs $(A',I')$ and $(A'',I'')$ on $X$
are equivalent closing subalgebras of $\mathcal{O}(X)$. 
Let $A''' \subset \mathcal{O}(X)$ be any subalgebra of finite type 
containing $A' \cup A''$. 
Define an ideal in $A'''$ by $I''' := \sqrt{\bangle{I'}}$. 
Then Lemma~\ref{naturalpairs} tells us that
the pair $(A''', I''')$ is a closing subalgebra.

We have to show that $I'''$ equals $\sqrt{\bangle{I''}}$. 
By Lemma~\ref{naturalpairs}, the morphism $X \to X'''$
defined by the inclusion $A''' \subset \mathcal{O}(X)$ is an open
embedding and $I''' \subset A'''$ is the vanishing ideal of 
$X''' \setminus X$. For every $f \in I''$, the map $X \to X'''$
restricts to an isomorphism $X_{f} \to X'''_{f}$. Hence the desired
identity of ideals follows from
$$ X = \bigcup_{f \in I''} X_{f}. $$

Finally, we show that if a closing subalgebra $(A'',I'')$ 
is equivalent to a natural pair $(A',I')$,
then also $(A'',I'')$ is natural. Choose $(A''',I''')$
as in~\ref{closedsubalgebra}~(ii). 
By Lemma~\ref{naturalpairs}~(ii), the pair $(A''',I''')$ 
is natural. In particular, $X_{f}$ is affine for every $f \in I''$.
Moreover, $X$ is covered by these $X_{f}$,
because $I'''$ equals $\sqrt{\bangle{I''}}$. \endproof

We are ready for the main result of this section.
Given a quasiaffine variety $X$, 
we denote as before by $\mathfrak{O}(X)$ 
the collection of all natural pairs on $X$.
For a morphism $\varphi \colon X \to Y$ of varieties, we denote by
$\varphi^{*} \colon \mathcal{O}(Y) \to \mathcal{O}(X)$ the pullback 
of functions. 

\begin{prop}\label{quasiaffequiv}
The assignments $X \mapsto (\mathcal{O}(X), \mathfrak{O}(X))$ and
$\varphi \mapsto \varphi^{*}$ define a contravariant equivalence of
the category of quasiaffine varieties with the category of
quasiaffine algebras. 
\end{prop}

\proof First of all, we check that the above assignment is
in fact well defined on morphisms. Let $\varphi \colon X \to Y$ be
any morphism of quasiaffine varieties. Choose a closing
subalgebra $(B',J')$ in $\mathfrak{O}(Y)$. 
By Lemma~\ref{naturalpairs}~(ii), we can construct a closing 
subalgebra $(A',I')$ in $\mathfrak{O}(X)$ such
that $\varphi^{*}(B') \subset A'$.

Now, consider the affine closures $X' := \Spec(A')$ and $Y' :=
\Spec(B')$ of $X$ and~$Y$. The morphism $\varphi' \colon X' \to 
Y'$ defined by the restriction $\varphi^{*} \colon B' \to A'$ maps
$X$ to~$Y$. Since $I'$ and $J'$ are precisely the vanishing ideals of
the complements $X' \setminus X$ and $Y' \setminus Y$, we obtain the
condition required in~\ref{quasiaffalgdef}~(ii):
$$I' \subset \sqrt{\bangle{\varphi^{*}(J')}}. $$  

Thus $\varphi \mapsto \varphi^{*}$ is in fact well defined. 
Moreover, $X \mapsto (\mathcal{O}(X), \mathfrak{O}(X))$ 
and $\varphi \mapsto \varphi^{*}$ clearly define a contravariant
functor, and this functor is injective on morphisms. 

For surjectivity, let $\mu \colon \mathcal{O}(Y) \to \mathcal{O}(X)$
be a homomorphism of quasiaffine algebras.
Let $(A',I') \in \mathfrak{O}(X)$ and $(B',J') \in \mathfrak{O}(Y)$ 
as in~\ref{quasiaffalgdef}~(ii).
Then $\mu$ defines a morphism $\varphi'$ from
$\Spec(A')$ to $\Spec(B')$.
The condition on the ideals and Lemma~\ref{naturalpairs}~(i)
ensure that $\varphi'$ restricts to 
a morphism $\varphi \colon X \to Y$. 
Clearly, we have $\varphi^{*} = \mu$. 

It remains to show that up to isomorphism, every quasiaffine algebra
$(A,\mathfrak{A})$ arises from a quasiaffine variety. 
Let $(A',I') \in \mathfrak{A}$, set $X' := \Spec(A')$, and let $X
\subset X'$ be the open subvariety obtained by removing the zero set
of $I'$. Then $\mathcal{O}(X) = A$, and $(A',I')$ is a natural pair on
$X$. Lemma~\ref{quasiaff2closingsubalg} gives $\mathfrak{O}(X) = 
\mathfrak{A}$. \endproof

We conclude this section with the observation, that restricted on
the category of quasiaffine varieties $X$ with $\mathcal{O}(X)$ of
finite type, our algebraic description collapses in a very convenient
way:

\begin{rem}\label{collaps}
For any quasiaffine algebra $(A, \mathfrak{A})$ we have
\begin{enumerate}
\item The algebra $A$ is of finite type over $\KK$ if and only if
  $(A,I) \in \mathfrak{A}$ holds with some radical ideal $I \subset
  A$.
\item The quasiaffine algebra $(A,\mathfrak{A})$ arises from an affine
  variety if and only if  $(A,A) \in \mathfrak{A}$ holds.
\end{enumerate}
\end{rem}

\section{Freely graded quasiaffine algebras}\label{section2}

In this section, we introduce the formal framework of homogeneous
coordinate rings, namely freely graded quasiaffine algebras and their
morphisms. The geometric interpretation of these notions amounts to an
equivariant version of the equivalence of categories presented in the
preceding section. 

\begin{defi}\label{freegradalgdef}
Let $(A,\mathfrak{A})$ be a quasiaffine algebra, and let $\Lambda$ be
a finitely generated abelian group. 
We say that $(A,\mathfrak{A})$ is {\em freely graded\/} by
$\Lambda$ (or {\em freely $\Lambda$-graded\/}) if there is 
a grading 
$$ A = \bigoplus_{L \in \Lambda} A_{L},  $$
and there exists a closing subalgebra $(A',I') \in \mathfrak{A}$
admitting homogeneous elements $f_{1}, \ldots, f_{r} \in I'$ such that 
$I'$ equals $\sqrt{\bangle{f_{1}, \ldots, f_{r}}}$
and every localization $A_{f_{i}}$ has in each degree $L \in
\Lambda$ a homogeneous invertible element. 
\end{defi}

\begin{exam}\label{polring}
For $n \ge 2$, the polynomial ring $\KK[T_{1}, \ldots, T_{n}]$
together with the usual $\ZZ$-grading can be made into a freely 
graded quasiaffine algebra: 
Let $\mathfrak{A}$ be the class of $(A,I)$, where
$I := \bangle{T_{1}, \ldots, T_{n}}$.
\end{exam}

The {\em weight monoid\/} of an integral domain $A$ graded by a
finitely generated abelian group $\Lambda$ is the submonoid
$\Lambda^{*} \subset \Lambda$ consisting of all weights $L \in
\Lambda$ with $A_{L} \ne \{0\}$. For the weight monoid of a freely
graded quasiaffine algebra, we have: 

\begin{rem}\label{pointedweightcone}
Let $(A,\mathfrak{A})$ be a freely $\Lambda$-graded quasiaffine
algebra. Then the weight monoid $\Lambda^{*} \subset \Lambda$
of $A$ generates $\Lambda$ as a group.
\end{rem}

We turn to homomorphisms. The final notion of a morphism of
freely graded quasiaffine algebras will be given below. 
First we have to consider homomorphisms that are compatible with 
the structure:  

\begin{defi}\label{admissiblehom}
Let the quasiaffine algebras $(A,\mathfrak{A})$ and $(B,\mathfrak{B})$
be freely graded by $\Lambda$ and $\Gamma$, respectively. A homomorphism
$\mu \colon (B,\mathfrak{B}) \to (A,\mathfrak{A})$ of quasiaffine
algebras is called {\em graded\/}, if there is a homomorphism $\t{\mu}
\colon \Gamma \to \Lambda$ with 
\begin{equation}
\label{gradedhomcond}
\mu(B_{E}) 
 \subset 
A_{\t{\mu}(E)} 
\quad \text{for all } E \in \Gamma.
\end{equation}  
\end{defi}

By Remark~\ref{pointedweightcone}, a graded homomorphism 
$\mu \colon (B,\mathfrak{B}) \to (A,\mathfrak{A})$ of freely graded
quasiaffine algebras uniquely determines its
accompanying homomorphism $\t{\mu} \colon \Gamma \to
\Lambda$. Moreover, the composition of two graded homomorphisms is
again graded.

For the subsequent treatment of our homogeneous coordinate rings we
need a coarser concept of a morphism of freely graded quasiaffine
algebras than the notion of a graded homomorphism would yield.
This is the following:

\begin{defi}\label{pointedmorphdef}
Let the quasiaffine algebras $(A,\mathfrak{A})$ and $(B,\mathfrak{B})$
be freely graded by finitely generated abelian groups $\Lambda$ and
$\Gamma$ respectively. 
\begin{enumerate}
\item Two graded homomorphisms $\mu, \nu \colon (B,\mathfrak{B})
  \to (A,\mathfrak{A})$ are called {\em equivalent\/} if there is a
  homomorphism $c \colon \Gamma \to A_{0}^{*}$ such that for every $E
  \in \Gamma$ and every $g \in B_{E}$ we have 
  $$ \nu(g) = c(E) \mu(g). $$
\item A {\em morphism\/} $(B,\mathfrak{B}) \to (A,\mathfrak{A})$ of
  the freely graded quasiaffine algebras $(B,\mathfrak{B})$ and
  $(A,\mathfrak{A})$ is the equivalence class $[\mu]$ of a graded
  homomorphism $\mu \colon (B,\mathfrak{B}) \to (A,\mathfrak{A})$.
\end{enumerate}
\end{defi}

In the setting of~(i) we shall say that $\mu$ and $\nu$ 
{\em differ by a character\/} $c \colon \Gamma \to A_{0}^{*}$. 
Since equivalence of graded homomorphisms is compatible with
composition, this definition makes the freely graded quasiaffine
algebras into a category.

We give now a geometric interpretation of the above notions. 
We assume for the rest of this section that if $\KK$ 
is of characteristic $p > 0$,
then our finitely generated abelian groups $\Lambda$ 
have no $p$-torsion, i.e.~$\Lambda$ contains no elements of order $p$. 
Under this assumption, each $\Lambda$ defines a
diagonalizable algebraic group
$$ H := \Spec(\KK[\Lambda]). $$

Recall that the characters of this group $H$ are precisely the
canonical generators $\chi^{L}$, $L \in \Lambda$, of the group
algebra $\KK[\Lambda]$. 
In fact, the assignment $\Lambda \mapsto H$ defines a
contravariant equivalence of categories, see for
example~\cite[Section~III.~8]{Bor}. 

Now, suppose that a diagonalizable group $H = \Spec(\KK[\Lambda])$
acts by means of a regular map $H \times X \to X$ on a (not
necessarily affine) variety $X$. A function $f \in \mathcal{O}(X)$ is
called {\em homogeneous\/} with respect to a character 
$\chi^{L} \colon H \to   \KK^{*}$ if for every $(t,x) \in H \times X$
we have
$$ f(t \mal x) = \chi^{L}(t) f(x). $$

For $L \in \Lambda$, let $\mathcal{O}(X)_{L} \subset \mathcal{O}(X)$
denote the subset of all $\chi^{L}$-homogeneous functions. 
It is well known, use for example \cite[p.~67~Lemma]{Kn}, that the
action of $H$ on $X$ defines a grading
$$ \mathcal{O}(X) = \bigoplus_{L \in \Lambda} \mathcal{O}(X)_{L}.$$

Recall that one obtains in this way a canonical correspondence
between affine $H$-varieties and $\Lambda$-graded affine algebras
(the arguments presented in~\cite[p.~11]{Do} for the case
$\Lambda = \ZZ$ also work in the general case).

We are interested in free $H$-actions on quasiaffine varieties $X$, 
where {\em free\/} means that every orbit map $H \to H \mal x$
is an isomorphism. In this situation, we have:

\begin{lemma}\label{gradedqavar2qaalg}
Let the diagonalizable group $H = \Spec(\KK[\Lambda])$ act 
freely by means of a regular map $H \times X \to X$ on a quasiaffine
variety $X$. 
Then the associated $\Lambda$-grading of $\mathcal{O}(X)$ makes
$(\mathcal{O}(X),\mathfrak{O}(X))$ into a freely graded quasiaffine
algebra.
\end{lemma}

\proof Let $(A'',I'')$ be any natural pair on $X$, and let $g_{1},
\ldots, g_{s}$ be a system of generators of $A''$. 
Let $A' \subset \mathcal{O}(X)$ denote the subalgebra generated 
by all the homogeneous components of the $g_{j}$. 
Then $A'$ is graded, and according to Lemma~\ref{naturalpairs}~(ii), 
we obtain a natural pair $(A',I')$ on $X$ by defining 
$I' := \sqrt{\bangle{I''}}$.

Now, the $\Lambda$-grading of $A'$ comes from an $H$-action on 
$X' := \Spec(A')$. This $H$-action extends the initial $H$-action on
$X$. In particular, the ideal $I' \subset A'$ is graded, because it is
the vanishing ideal of the invariant set $X' \setminus X$, see
Lemma~\ref{naturalpairs}~(i). This fact enables us to verify
the condition of~\ref{freegradalgdef} for  $I'$:

\goodbreak

Choose generators $L_{1}, \ldots, L_{k}$ of $\Lambda$. 
Consider $x \in X$, and choose a homogeneous 
$h \in I'$ with $h(x) \ne 0$.
Since $H$ acts freely, the orbit map $H \to H \mal x$
is an isomorphism. 
Thus we find for every $i$ a $\chi^{L_{i}}$-homogeneous 
regular function $h_{i}$ on  $H \mal x$ with $h_{i}(x) \ne 0$. 
Since $H \mal x$ is closed in $X_{h}$, the
$h_{i}$ extend to $\chi^{L_{i}}$-homogeneous regular
functions on $X_{h}$. 

For a suitable $r > 0$, the product $f := h^{r}h_{1} \ldots h_{k}$
is a regular function on $X'$ with $f \in \bangle{h}$ and hence
$f \in I'$. 
By construction, $f$ is homogeneous, and we have $f(x) \ne 0$. 
Moreover, the Laurent monomials
in $h_{1}, \ldots, h_{k}$ provide for each degree $L \in \Lambda$
a $\chi^{L}$-homogeneous invertible function on $X_{f}$. 
Since finitely many of the $X_{f}$ cover $X$, this gives the desired
property on the ideal $I' \subset A'$.  
\endproof

In order to give the equivariant version of
Proposition~\ref{quasiaffequiv}, we have to fix the notion of a
morphism of quasiaffine varieties with an action of a diagonalizable
group. This is the following: 

\begin{defi}
Let $G \times X \to X$ and $H \times Y \to Y$ be algebraic group
actions. A morphism $\varphi \colon X \to Y$ is called 
{\em equivariant\/} if there is a homomorphism $\t{\varphi} \colon G
\to H$ of algebraic groups such that for all
$(g,x) \in G \times X$ we have 
$$\varphi(g \mal x) = \t{\varphi}(g) \mal \varphi(x).$$
\end{defi}

This notion of an equivariant morphism makes the quasiaffine varieties
with a free action of a diagonalizable group into a category. We
obtain the following equivariant version of
Proposition~\ref{quasiaffequiv}: 

\begin{prop}\label{equivquasiaffequiv}
The assignments $X \mapsto (\mathcal{O}(X), \mathfrak{O}(X))$ and
$\varphi \mapsto \varphi^{*}$ define a contravariant equivalence from
the category of quasiaffine varieties with a free diagonalizable group
action to the category of freely graded quasiaffine algebras and
graded homomorphisms.
\end{prop}

\proof By Lemma~\ref{gradedqavar2qaalg}, the assignment $X \mapsto
(\mathcal{O}(X), \mathfrak{O}(X))$ is well defined. From 
Proposition~\ref{quasiaffequiv} and the observation that 
equivariant morphisms of quasiaffine varieties correspond to graded 
homomorphisms of quasiaffine algebras we infer functoriality and
bijectivity on the level of morphisms.

In order to see that up to isomorphism any quasiaffine
algebra $(A,\mathfrak{A})$ which is freely graded by some $\Lambda$
arises in the above manner from a quasiaffine variety with
free diagonalizable group action, we repeat the corresponding part of
the proof of Proposition~\ref{quasiaffequiv} in an equivariant manner:

Let $(A',I')$ be as in Definition~\ref{freegradalgdef}. Let 
$A'' \subset A$ be any graded subalgebra of finite type with 
$A' \subset A''$, and let $I'' := \sqrt{\bangle{I'}}$.
Then $(A'',I'')$ belongs to $\mathfrak{A}$, and 
the ideal $I''$ still satisfies 
the condition of Definition~\ref{freegradalgdef}. 

The affine variety $X'' := \Spec(A'')$ comes along with an action of
the diagonalizable group $H := \Spec(\KK[\Lambda])$ such that the
corresponding grading of $\mathcal{O}(X'') = A''$ gives back the
original $\Lambda$-grading of the algebra $A''$. 
Removing the $H$-invariant zero set of $I''$ from $X''$, 
gives a quasiaffine $H$-variety $X$. 

By construction, the $\Lambda$-graded algebras $\mathcal{O}(X)$ and
$A$ coincide, and $(A'',I'')$ is a natural pair on $X$. Moreover, the
local existence of invertible homogeneous functions in each degree
implies that for every $x \in X$ the orbit map $H \mapsto H \mal x$ is
an isomorphism. In other words, the action of $H$ on $X$ is free. 
\endproof

\begin{exam}
The standard $\KK^{*}$-action on $\KK^{n+1} \setminus \{0\}$ has
$(A,\mathfrak{A})$ of Example~\ref{polring} as associated freely
graded quasiaffine algebra. 
\end{exam}

The remaining task is to translate the notion of equivalence of
graded homomorphisms. For this let $X$ and $Y$ be quasiaffine 
varieties with actions of diagonalizable groups 
$H := \Spec(\KK[\Lambda])$ and $G := \Spec(\KK[\Gamma])$. 
Denote by $(A,\mathfrak{A})$ and $(B,\mathfrak{B})$ the freely 
graded quasiaffine algebras associated to $X$ and $Y$. 

\begin{rem}\label{equivgeom}
Two graded homomorphisms $\mu, \nu \colon (B,\mathfrak{B}) \to
(A,\mathfrak{A})$ are equivalent
if and only if there is an $H$-invariant
morphism $\gamma \colon X \to G$ such that
the morphisms $\varphi, \psi \colon X \to Y$ corresponding to $\mu$
and $\nu$ always satisfy $\psi(x) = \gamma(x) \mal \varphi(x)$.
\end{rem}

\section{Picard graded sheaves of algebras}\label{section3}

Let $X$ be an algebraic variety and denote by $\Pic(X)$
its Picard group.
In this section we prepare the definition of a 
graded ring structure on the vector space
$$ \bigoplus_{[L] \in \Pic(X)} H^{0}(X,L). $$

More generally, we even need a ring structure for the 
corresponding sheaves of vector spaces.
The problem is easy, if $\Pic(X)$ is free: 
Then we can realize it as a group $\Lambda$ 
of line bundles as in~\cite[Sec.~2]{Ha},
and we can work with the associated
$\Lambda$-graded $\mathcal{O}_{X}$-algebra $\mathcal{R}$.

If $\Pic(X)$ has torsion, then we can at most expect
a surjection $\Lambda \to \Pic(X)$ with a free group
$\Lambda$ of line bundles.
Thus the problem is to identify in a suitable manner
isomorphic homogeneous
components of the $\Lambda$-graded $\mathcal{O}_{X}$-algebra 
$\mathcal{R}$. 
This is done by means of shifting families
and their associated ideals
$\mathcal{I} \subset \mathcal{R}$, see~\ref{shiftfamdef}
and~\ref{associdealdef}.
The quotient $\mathcal{A} := \mathcal{R} / \mathcal{I}$
then will realize the desired ring structure.

To begin, let us recall the necessary constructions 
from~\cite{Ha}. 
Consider an open cover $\mathfrak{U} = (U_{i})_{i \in I}$ of~$X$. 
This cover gives rise to an additive group $\Lambda(\mathfrak{U})$ 
of line bundles on~$X$: For each \v{C}ech cocycle 
$\xi \in Z^{1}(\mathfrak{U},\mathcal{O}_{X}^{*})$, let
$L_{\xi}$ denote the line bundle obtained by gluing the products
$U_{i} \times \KK$ along the maps
$$ (x,z) \mapsto (x, \xi_{ij}(x)z).$$
Define the sum $L_{\xi} + L_{\eta}$ of two such line bundles to be
$L_{\xi\eta} = L_{\eta\xi}$. This makes the set $\Lambda(\mathfrak{U})$
consisting of all the bundles $L_{\xi}$ into an abelian group,
which is isomorphic to the cocycle group 
$Z^{1}(\mathfrak{U},\mathcal{O}_{X}^{*})$. 

When we speak of a {\em group of line bundles\/} on $X$, 
we think of a finitely
generated free subgroup of some group $\Lambda(\mathfrak{U})$ as
above.
Note that for any such group $\Lambda$ of line bundles, we have a
canonical homomorphism $\Lambda \to \Pic(X)$ to the Picard group.

We come to the graded $\mathcal{O}_{X}$-algebra associated to a
group $\Lambda$ of line bundles on a variety $X$. For each line
bundle $L \in \Lambda$, let $\mathcal{R}_{L}$ denote its sheaf of
sections. In the sequel, we shall identify $\mathcal{R}_{0}$ with the
structure sheaf $\mathcal{O}_{X}$. The 
{\em graded $\mathcal{O}_{X}$-algebra\/} associated to $\Lambda$ is
the quasicoherent sheaf
$$ \mathcal{R} := \bigoplus_{L \in \Lambda} \mathcal{R}_{L}, $$
where the multiplication is defined as follows: The sections of a
bundle $L_{\xi} \in \Lambda$ over an open set 
$U \subset X$ are described by families $f_{i} \in \mathcal{O}_{X}(U
\cap U_{i})$ that are compatible with the cocycle $\xi$. For
any two sections $f \in \mathcal{R}_{L}(U)$ and $f' \in
\mathcal{R}_{L'}(U)$, the product $(f_{i}f'_{i})$ of their
defining families $(f_{i})$ and $(f'_{i})$ gives us a
section $ff' \in \mathcal{R}_{L+L'}(U)$. 

In the sequel, we fix an open cover 
$\mathfrak{U} = (U_{i})_{i \in I}$ of~$X$
and a group $\Lambda \subset \Lambda(\mathfrak{U})$ 
of line bundles. 
Let $\mathcal{R}$ denote the associated $\Lambda$-graded
$\mathcal{O}_{X}$-algebra.
Here comes the notion of a shifting family for $\mathcal{R}$:

\begin{defi}\label{shiftfamdef}
Let $\Lambda_{0} \subset \Lambda$ be any subgroup of the kernel 
of $\Lambda \to \Pic(X)$. 
By a {\em $\Lambda_{0}$-shifting family\/} for $\mathcal{R}$ we
mean a family $\varrho = ( \varrho_{E} )$ of $\mathcal{O}_{X}$-module
isomorphisms $\varrho_{E} \colon \mathcal{R} \to \mathcal{R}$, where $E \in
\Lambda_{0}$, with the following properties:
\begin{enumerate}
\item for every $L \in \Lambda$ and every $E \in \Lambda_{0}$ the
  isomorphism $\varrho_{E}$ maps $\mathcal{R}_{L}$ onto
  $\mathcal{R}_{L+E}$, 
\item for any two $E_{1}, E_{2} \in \Lambda_{0}$ we have
  $\varrho_{E_{1} + E_{2}} = \varrho_{E_{2}} \circ \varrho_{E_{1}}$,
\item for any two homogeneous sections $f,g$ of $\mathcal{R}$ and
  every $E \in \Lambda_{0}$ we have 
  $\varrho_{E}(fg) = f \varrho_{E}(g)$. 
\end{enumerate}
If $\Lambda_{0}$ is the full kernel of $\Lambda \to \Pic(X)$, then we
also speak of a {\em full shifting family\/} for $\mathcal{R}$ instead
of a $\Lambda_{0}$-shifting family.
\end{defi}

The first basic observation is existence of shifting families and a
certain uniqueness statement: 

\begin{lemma}\label{shiftfamprops}
Let $\Lambda_{0} \subset \Lambda$ be a subgroup of the kernel of 
$\Lambda \to \Pic(X)$. Then there exist $\Lambda_{0}$-shifting
families for $\mathcal{R}$, and any two such families $\varrho$,
$\varrho'$ differ by a character $c \colon \Lambda_{0} \to
\mathcal{O}^{*}(X)$ in the sense that $\varrho'_{E} = c(E)\varrho_{E}$
holds for all $E \in \Lambda_{0}$.  
\end{lemma}

\proof 
For the existence statement, fix a $\ZZ$-basis of the subgroup 
$\Lambda_{0} \subset \Lambda$. 
For any member $E$ of this basis choose a bundle isomorphism 
$\alpha_{E} \colon 0 \to E$ from the trivial bundle $0 \in \Lambda$
onto $E \in \Lambda$. 
With respect to the cover $\mathcal{U}$, this
isomorphism is fibrewise multiplication with certain 
$\alpha_{i} \in \mathcal{O}^{*}(U_{i})$; so, on 
$U_{i} \times \KK$ it is of the form
\begin{equation}\label{localdata}
(x,z) \mapsto (x, \alpha_{i}(x)z).
\end{equation}

If $\alpha_{E'} \colon 0 \to E'$ denotes the isomorphism for a
further member of the basis of $\Lambda_{0}$, then the products
$\alpha_{i}\alpha_{i}'$  
of the corresponding local data define an isomorphism 
$\alpha_{E+E'} \colon 0 \to E + E'$. 
Similarly, by inverting local data, we obtain isomorphisms
$\alpha_{-E} \colon 0 \to -E$. 
Proceeding this way, we obtain an isomorphism
$\alpha_{E} \colon 0 \to E$ for every $E \in \Lambda_{0}$.

The local data $\alpha_{i}$ of an isomorphism 
$\alpha_{E} \colon 0 \to E$ as constructed above define as well an
isomorphism $L \to L + E$ for any $L \in \Lambda$. 
By shifting homogeneous sections according to 
$f \mapsto \alpha_{E} \circ f$,
one obtains $\mathcal{O}_{X}$-module 
isomorphisms $\varrho_{E} \colon \mathcal{R} \to \mathcal{R}$ mapping
each $\mathcal{R}_{L}$ onto $\mathcal{R}_{L+E}$. The
Properties~\ref{shiftfamdef}~(ii) and~(iii) are then clear by
construction. 

We turn to the uniqueness statement. Let $\varrho$, $\varrho'$ be two
$\Lambda_{0}$-shifting families for $\mathcal{R}$.
Using Property~\ref{shiftfamdef}~(iii) we
see that for every $E \in \Lambda_{0}$ and every homogeneous section
$f$ of $\mathcal{R}$, we have
$$ \varrho_{E}^{-1} \circ \varrho'_{E} (f) = \varrho_{E}^{-1} \circ
\varrho'_{E} (f \cdot 1) = f \cdot \varrho_{E}^{-1} \circ \varrho'_{E} (1). $$
Thus, setting $c(E) := \varrho_{E}^{-1} \circ \varrho'_{E} (1)$ we obtain a
map $c \colon \Lambda_{0} \to \Of^{*} (X)$ such that $\varrho'_{E}$ equals
$c(E)\varrho_{E}$. Properties~\ref{shiftfamdef}~(ii) and~(iii) show
that $c$ is a homomorphism:
\begin{eqnarray*}
c(E_{1}+E_{2}) 
& = & 
\varrho_{E_{1}+E_{2}}^{-1} \circ \varrho'_{E_{1}+E_{2}}(1) \\
& = & 
\varrho_{-E_{1}} \circ \varrho_{-E_{2}} \circ  \varrho'_{E_{2}} 
\circ \varrho'_{E_{1}}(1) \\
& = & 
\varrho_{-E_{1}} \circ \varrho_{-E_{2}} \circ \varrho'_{E_{2}} 
(\varrho'_{E_{1}}(1) \mal 1)\\
& = & 
\varrho_{-E_{1}} ( \varrho'_{E_{1}}(1) c(E_{2})) \\
& = & c(E_{1})c(E_{2}). \qquad \qed
\end{eqnarray*}

\goodbreak

We shall now associate to any shifting family an ideal in the
$\mathcal{O}_{X}$-algebra $\mathcal{R}$. First we remark
that for any subgroup $\Lambda_{0} \subset \Lambda$ the algebra
$\mathcal{R}$ becomes $\Lambda/\Lambda_{0}$-graded by defining 
the homogeneous component of a class $[L] \in \Lambda/\Lambda_{0}$
as 
$$ \mathcal{R}_{[L]} := \sum_{L' \in [L]} \mathcal{R}_{L'}.$$

\begin{lemma}\label{associdealprops}
Let $\Lambda_{0}$ be a subgroup of the kernel of $\Lambda \to
\Pic(X)$, and let $\varrho$ be a $\Lambda_{0}$-shifting family.
For each given open subset $U \subset X$ consider the ideal
$$ 
\mathfrak{I} (U) 
\; := \; 
\bangle{f - \varrho_{E}(f); \; 
f \in \mathcal{R}(U), \;
E \in \Lambda_{0}} 
\; \subset \; 
\mathcal{R}(U).
$$
Let $\mathcal{I}$ denote the sheaf associated to
the presheaf $U \mapsto \mathfrak{I}(U)$. Then $\mathcal{I}$ is
a quasicoherent ideal of $\mathcal{R}$, and we have:
\begin{enumerate}
\item Every $\mathcal{I}(U)$ is homogeneous with respect to the
$\Lambda/\Lambda_{0}$-grading of $\mathcal{R}(U)$. 
\item For every $L \in \Lambda$ we have 
$ \mathcal{R}_{L}(U) \cap \mathcal{I}(U) = \{0\}$.
\end{enumerate}
\end{lemma}

\proof 
First note that the ideal sheaf $\mathcal{I}$ is indeed quasicoherent,
because it is a sum of images of quasicoherent sheaves.

We check~(i). 
Using Property~\ref{shiftfamdef}~(iii), we see that 
each ideal $\mathfrak{I}(U)$ is generated by the elements 
$1 - \varrho_{E}(1)$,
where $E \in \Lambda_{0}$. 
Consequently, each stalk $\mathfrak{I}_{x}$ is a 
$\Lambda/\Lambda_{0}$-homogeneous ideal in $\mathcal{A}_{x}$.
This implies that the associated sheaf $\mathcal{I}$ is a 
$\Lambda/\Lambda_{0}$-homogeneous ideal sheaf in  $\mathcal{A}$.

We turn to~(ii). By construction, it suffices to
consider local sections $f \in \mathcal{R}_{L}(U)
\cap \mathfrak{I}(U)$. By the definition of $\mathfrak{I}(U)$ and
Property~\ref{shiftfamdef}~(iii), there exist homogeneous elements
$f_{i} \in \mathcal{R}_{L_{i}}(U)$ such that we can write $f$ as 
\begin{equation}\label{minimalrep}
f = \sum_{i=1}^{r} f_{i} - \varrho_{E_{i}}(f_{i}).
\end{equation}

Since $\mathfrak{I}(U)$ is $\Lambda/\Lambda_{0}$-graded, 
all the $L_{i}$ belong to the class $[L]$ in
$\Lambda/\Lambda_{0}$. Moreover, we can achieve in the
representation~(\ref{minimalrep}) of $f$ that all $f_{i}$ 
are of degree $L \in [L]$. 
Namely, we can use Property~\ref{shiftfamdef}~(ii) to write 
$f_{i} - \rho_{E_{i}} (f_{i})$ in the form
\begin{eqnarray*}
f_{i} - \varrho_{E_{i}} (f_{i}) 
&  = &  
\varrho_{L-L_{i}} (f_{i}) - 
\varrho_{E_{i} + L_{i} - L} (\varrho_{L-L_{i}} (f_{i})) 
\\  
& & + (- \varrho_{L-L_{i}} (f_{i})) -
\varrho_{L_{i}-L} (-\varrho_{L - L_{i}}(f_{i})).
\end{eqnarray*}

Moreover we can choose the representation~(\ref{minimalrep}) minimal
in the sense that $r$ is minimal with the property that every $f_{i}$
is of degree $L$. Then the $E_{i}$ are pairwise different from each
other, because otherwise we could shorten the representation by
gathering. But this implies $\varrho_{E_{i}}(f_{i}) = 0$ for every
$i$. Hence we obtain $f=0$. 
\endproof

\begin{defi}\label{associdealdef}
Let $\Lambda_{0}$ be a subgroup of the kernel of $\Lambda \to
\Pic(X)$, and let $\varrho$ be a $\Lambda_{0}$-shifting family
for $\mathcal{R}$.
The {\em ideal associated to $\varrho$} is the
$\Lambda/\Lambda_{0}$-graded ideal sheaf $\mathcal{I}$ of
$\mathcal{R}$ defined in~\ref{associdealprops}.
\end{defi}

With the aid of the ideal associated to a shifting family, we can pass from
$\mathcal{R}$ to more coarsely graded $\mathcal{O}_X$-algebras:

\begin{lemma}\label{gradproject}
Let $\Lambda_0 \subset \Lambda$ be a subgroup, and let $\varrho$ be a
$\Lambda_0$-shifting family with associated ideal $\mathcal{I}$. Set
$\mathcal{A} := \mathcal{R}/\mathcal{I}$, and let $\pi \colon
\mathcal{R} \to \mathcal{A}$ denote the projection.
\begin{enumerate}
\item The $\mathcal{O}_{X}$-algebra $\mathcal{A}$ is quasicoherent,
  and it inherits a
  $\Lambda/\Lambda_{0}$-grading from $\mathcal{R}$ as follows
 $$ \mathcal{A} 
    = \bigoplus_{[L] \in \Lambda/\Lambda_{0}} \mathcal{A}_{[L]}
    := \bigoplus_{[L] \in \Lambda/\Lambda_{0}}
    \pi(\mathcal{R}_{[L]}). $$
\item For any $L \in \Lambda$ the induced map $\pi_{L} \colon \mathcal{R}_{L}
\to \mathcal{A}_{[L]}$ is an isomorphism of
$\mathcal{O}_{X}$-modules. In particular, we obtain
$$\mathcal{A}(X) \cong \mathcal{R}(X) / \mathcal{I}(X).$$
\item The $\mathcal{O}_{X}$-algebra $\mathcal{A}$ is locally generated
  by finitely many invertible homogeneous elements.
\end{enumerate}
\end{lemma}

\proof 
The first assertion follows directly from the fact that we have a
commutative diagram where the lower arrow is an isomorphism of sheaves:
$$ 
\xymatrix{
 & {\mathcal{R}} \ar[ld]_{\pi} \ar[rd] & \\
{\mathcal{A}} \ar[rr] & & 
{\bigoplus_{[L] \in \Lambda/\Lambda_{0}} 
 \mathcal{R}_{[L]}/\mathcal{I}_{[L]}} 
}
$$

To prove (ii), note that $\pi_{L} \colon \mathcal{R}_{L} \to
\mathcal{A}_{[L]}$ is injective by Lemma~\ref{associdealprops}~(ii).
For bijectivity, we have to show that $\pi_{L}$ is stalkwise surjective.
Let $h$ be a local section of ${\mathcal{A}_{[L]}}$ near some $x \in X$. 
Since ${\mathcal{A}_{[L]}}$ equals $\pi(\mathcal{R}_{[L]})$, 
we may assume that $h = \pi(f)$ with a local section $f$ of 
$\mathcal{R}_{[L]}$ near $x$. Write $f$ as the sum of its
$\Lambda$-homogeneous components: 
$$ f = \sum_{L' \in [L]} f_{L'}. $$

For every $L' \ne L$, we subtract $f_{L'} - \varrho_{L-L'}(f_{L'})$
from $f$. The result is a local section $g$ of $\mathcal{R}_{L}$ near
$x$ which still projects onto $h$. This proves bijectivity of 
$\pi_{L} \colon \mathcal{R}_{L} \to \mathcal{A}_{[L]}$. 
The isomorphy on the level of global sections then is due
to left exactness of the section functor.

To prove assertion~(iii), note that the analogous statement
holds for $\mathcal{R}$. 
In fact, for small $U \subset X$, the algebra
$\mathcal{R}(U)$ is even a Laurent monomial algebra over
$\mathcal{O}(U)$. 
Together with assertion (ii), this observation
gives statement~(iii).
\endproof

\begin{defi}\label{picgradalg}
Let $\Lambda_{0} \subset \Lambda$ be a subgroup of the kernel $\Lambda \to
\Pic(X)$, and let $\varrho$ be a $\Lambda_{0}$-shifting family for
$\mathcal{R}$ with associated ideal $\mathcal{I}$.  We call the
$\Lambda/\Lambda_{0}$-graded $\mathcal{O}_{X}$-algebra $\mathcal{A} :=
\mathcal{R}/ \mathcal{I}$ of~\ref{gradproject} the 
{\em Picard graded algebra\/} associated to $\varrho$.
\end{defi}

If every global invertible function on $X$ is constant, 
then the Picard graded algebras associated to different
$\Lambda_{0}$-shifting families are isomorphic 
(a graded homomorphism of sheaves is defined by 
requiring~\ref{gradedhomcond} on the level of sections):

\begin{lemma}\label{shiftfamunique}
Suppose $\mathcal{O}^{*}(X) = \KK^{*}$.
Let $\Lambda_{0} \subset \Lambda$ be a subgroup,
and let $\varrho$, $\varrho'$ be $\Lambda_{0}$-shifting 
families for $\mathcal{R}$ with associated ideals $\mathcal{I}$ and 
$\mathcal{I}'$. Then there is a graded automorphism of 
$\mathcal{R}$ having the identity of $\Lambda$ as accompanying
homomorphism and mapping $\mathcal{I}$ onto $\mathcal{I}'$.
\end{lemma}

\proof 
By Lemma~\ref{shiftfamprops}, there exists a
homomorphism $c \colon E \to \KK^{*}$ such that
$\varrho'_{E} = c(E) \varrho_{E}$ holds. By Lemma~\ref{charext} stated
below, this homomorphism extends to a homomorphism 
$c \colon \Lambda \to \KK^{*}$. Thus we can define the
desired automorphism $\mathcal{R} \to \mathcal{R}$ by mapping a 
section $f \in \mathcal{R}_{L}(U)$ to $c(L)f \in \mathcal{R}_{L}(U)$. 
\endproof

In the proof of this lemma, we made use of the following
standard property of lattices:

\begin{lemma}\label{charext}
Let $\Lambda_{0} \subset \Lambda$ be an inclusion of lattices. Then
any homomorphism $\Lambda_{0} \to \KK^{*}$ extends to a homomorphism 
$\Lambda \to \KK^{*}$.
\end{lemma}

Let us give a geometric interpretation of Picard graded algebras. 
Let $\Lambda$ be a group of
line bundles on $X$ with associated $\Lambda$-graded
$\mathcal{O}_{X}$-algebra $\mathcal{R}$. Fix a subgroup $\Lambda_{0}$
of the kernel of $\Lambda \to \Pic(X)$ and a $\Lambda_{0}$-shifting
family $\varrho$ for $\mathcal{R}$. 

Similar to the preceding section, we assume for the rest of this
section that in the case of a ground field $\KK$ of characteristic 
$p> 0$, the group $\Lambda/\Lambda_{0}$ has no $p$-torsion.
Under this hypothesis, we can show that the quotient 
$\mathcal{A} := \mathcal{R}/\mathcal{I}$ by the ideal 
associated to the shifting family $\varrho$ is reduced:

\begin{lemma}\label{reduced}
For every open $U \subset X$, the ideal $\mathcal{I}(U)$ is a radical
ideal in $\mathcal{R}(U)$. 
\end{lemma}

\proof First note that we may assume that $U$ is a small affine open
set such that $\mathcal{R}(U)$ is of finite type. Consider
the affine variety $Z := \Spec(\mathcal{R}(U))$. 
Then the $\Lambda/\Lambda_{0}$-grading of $\mathcal{R}(U) =
\mathcal{O}(Z)$ defines an action of the diagonalizable group 
$H := \Spec(\KK[\Lambda/\Lambda_{0}])$ on $Z$. Let 
$Z_{0} \subset Z$ denote the zero set of the ideal $\mathcal{I}(U)
\subset \mathcal{R}(U)$.

Now we can enter the proof of the assertion. 
Let $f \in \mathcal{O}(Z)$ with $f^{n}
\in \mathcal{I}(U)$. We have to show that $f \in \mathcal{I}(U)$
holds. Consider the decomposition of $f$ into homogeneous parts:
$$ f = \sum_{[L] \in \Lambda/\Lambda_{0}} f_{[L]}. $$
Since $f$ vanishes along the $H$-invariant zero set $Z_{0}$ of the
$\Lambda/\Lambda_{0}$-graded ideal $\mathcal{I}(U)$, also every
homogeneous component $f_{[L]}$ has to vanish along $Z_{0}$. 

We show that every $f_{[L]}$ belongs to $\mathcal{I}(U)$. Since
the $f_{[L]}$ vanish along $Z_{0}$, Hilbert's Nullstellensatz tells us
that for every degree $[L]$ some power $f_{[L]}^{m}$ lies in
$\mathcal{I}(U)$. Now consider
$$ g := \sum_{L' \in [L]} f_{L'} - (f_{L'} -
\varrho_{L-L'}(f_{L'})). $$

Then $g$ is $\Lambda$-homogeneous of degree $L$. Moreover, by explicit
multiplication, we see $g^{m} \in \mathcal{I}(U)$. But any
$\Lambda$-homogeneous element of $\mathcal{I}(U)$ is trivial. Thus
$g^{m} = 0$. Hence $g=0$, which in turn implies 
$f_{[L]} \in \mathcal{I}(U)$.
\endproof

In our geometric interpretation, we use the global
``$\Spec$''-construction, see for example~\cite{Ht}. 
Moreover, for any homogeneous
section $f \in \mathcal{A} (U)$, we denote 
its zero set in $X$ by $Z(f)$. This is well defined, because the
components $\mathcal{A}_{[L]}$ are locally free due to
Lemma~\ref{gradproject}~(ii). 

\begin{prop}\label{geominterp}
Let $\rq{X} := \Spec(\mathcal{A})$, and let $q \colon \rq{X} \to X$
be the canonical map.
\begin{enumerate}
\item $\rq{X}$ is a variety, $q \colon \rq{X} \to X$ is an
  affine morphism, and we have $\mathcal{A} = q_{*}
  \mathcal{O}_{\rq{X}}$.
\item For a homogeneous section $f \in \mathcal{A}_{[L]}(X)$ we obtain
  $q^{-1}(Z(f)) = V(\rq{X};f)$, where $V(\rq{X};f)$ is
  the zero set of the function  $f \in \mathcal{O}(\rq{X})$.
\item If $f_{i} \in \mathcal{A}(X)$ are homogeneous sections such that
  the sets $X \setminus Z(f_{i})$ are affine and cover $X$, then
  $\rq{X}$ is a quasiaffine variety.
\end{enumerate}
\end{prop}

\proof To check~(i), note that $\rq{X}$ is indeed a variety, because
by Lemmas~\ref{gradproject}~(iii) and~\ref{reduced}, the algebra
$\mathcal{A}$ is reduced and locally of finite type.
The rest of~(i) are standard properties of the
global ``$\Spec$''-construction for sheaves of
$\mathcal{O}_{X}$-algebras. 

Assertion~(ii) is clear in the case $[L] = 0$, because then we have
$\mathcal{A}_{0} = \mathcal{O}_{X}$. For a general $[L]$, we may
reduce to the previous case by multiplying $f$ locally with invertible
sections of degree $-[L]$. Note that invertible sections exist locally
by Lemma~\ref{gradproject}~(iii).
\endproof


\section{The homogeneous coordinate ring}\label{section4}

In this section, we give the precise definition of 
the homogeneous coordinate ring of a given variety, 
see Definition~\ref{homcoorddef}.
Moreover, we show in Proposition~\ref{uniquehomcoord}
that the homogeneous coordinate ring is
unique up to isomorphism.

In order to fix the setup,
recall from~\cite{Bo} that a (neither necessarily separated nor
irreducible) variety $X$, is said to be {\em divisorial\/} if every $x
\in X$ admits an affine neighbourhood of the form $X \setminus
Z(f)$ where $Z(f)$ is the zero set of a global section $f$ of some
line bundle $L$ on~$X$.  

\begin{rem}
Every separated irreducible $\QQ$-factorial variety is divisorial, and
every quasiprojective variety is divisorial.
\end{rem}

Here is the setup of this section: 
We assume that the multiplicative group $\KK^{*}$ 
is of infinite rank over $\ZZ$,
e.g.~$\KK$ is of characteristic zero or it is uncountable.
The variety $X$ is divisorial and satisfies
$\mathcal{O}^{*}(X) = \KK^{*}$.
Moreover, $\Pic(X)$ is finitely
generated and, if $\KK$ is of characteristic $p > 0$,
then $\Pic(X)$ has no $p$-torsion.

\begin{lemma}\label{ontopic}
There exists a group $\Lambda$ of line bundles on $X$ mapping
onto $\Pic(X)$. For any such $\Lambda$ the associated
$\Lambda$-graded $\mathcal{O}_{X}$-algebra $\mathcal{R}$ admits
homogeneous global sections $h_{1}, \ldots, h_{r}$ such that the sets 
$X \setminus Z(h_{i})$ are affine and cover~$X$.
\end{lemma} 

\proof
Only for the first statement there is something to show.
For this, we may assume that $\Pic(X)$ is not trivial.
Write $\Pic(X)$ as a direct sum of cyclic groups 
$\Pi_{1}, \ldots, \Pi_{m}$
and fix a generator $P_{l}$ for each $\Pi_{l}$. 
Choose a finite open cover $\mathfrak{U}$ of $X$ 
such that 
each $P_{l}$ is represented
by a cocycle 
$\xi^{(l)} \in Z^{1}(\mathfrak{U},\mathcal{O}^{*})$.
Choose members 
$U_{i}, U_{j}$ of $\mathfrak{U}$ such that 
$U_{i} \ne U_{j}$ holds and there is a point 
$x_{0} \in U_{i} \cap U_{j}$.

We adjust the $\xi^{(l)}$ as follows:
By the assumption on the ground field~$\KK$, we find
$a_{1}, \ldots, a_{m} \in \KK^{*}$ which are
linearly independent over $\ZZ$.
Define a locally constant cochain $\eta^{(l)}$ 
by setting $\eta^{(l)} := a_{l}/\xi^{(l)}_{ij}(x_{0})$ 
on $U_{i}$ and $\eta^{(l)} := 1$ on the $U_{k}$ different 
from $U_{i}$. 
Let $\zeta^{(l)} \in Z^{1}(\mathfrak{U},\mathcal{O}^{*})$ 
be the product of $\xi^{(l)}$ with the coboundary
of $\eta^{(l)}$.
 
Let $\Lambda \subset \Lambda(\mathfrak{U})$ be the 
subgroup generated by the line bundles arising from 
$\zeta^{(1)}, \ldots, \zeta^{(m)}$.
By construction $\Lambda$ maps onto $\Pic(X)$.
Moreover, we have 
$$ 
\Bigl(\bigl(\zeta_{ij}^{(1)}\bigr)^{n_{1}} \ldots
\bigl(\zeta_{ij}^{(m)}\bigr)^{n_{m}}\Bigr)(x_{o}) 
= a_{1}^{n_{1}} \ldots a_{m}^{n_{m}} $$
for the cocycle corresponding to a general element of 
$\Lambda$. By the choice of the $a_{l}$, this cocylce 
can only be trivial if all exponents $n_{l}$ vanish. 
It follows that $\Lambda$ is free.
\endproof

We fix a group $\Lambda$ of line bundles on $X$
as provided by Lemma~\ref{ontopic}, and a full shifting 
family $\varrho$ for the $\Lambda$-graded $\mathcal{O}_{X}$-algebra
$\mathcal{R}$ associated to $\Lambda$.
Let $\mathcal{I}$ denote
the ideal associated to the shifting family $\varrho$.
As seen in Lemma~\ref{gradproject}~(i), the $\mathcal{O}_{X}$-algebra
$\mathcal{A} := \mathcal{R}/\mathcal{I}$ is graded by $\Pic(X)$. 
In particular, we have a grading
$$ 
\mathcal{A}(X)
=
\bigoplus_{[L] \in \Pic(X)} \mathcal{A}_{[L]}(X).
$$

According to Lemmas~\ref{gradproject}~(ii) and~\ref{ontopic}, there 
are homogeneous $f_{1}, \ldots, f_{r} \in \mathcal{A}(X)$ such
that the sets $X \setminus Z(f_{i})$ are affine and cover $X$.
Hence Proposition~\ref{geominterp}~(iii) tells us that 
the variety $\rq{X} := \Spec(\mathcal{A})$ is quasiaffine.
Thus we have the collection $\mathfrak{A}(X)$ of 
natural pairs on $\rq{X}$ as closing subalgebras for 
$\mathcal{A}(X) = \mathcal{O}(\rq{X})$,   
see Lemma~\ref{quasiaff2closingsubalg}.

\begin{prop}\label{coringisqualg}
The pair $(\mathcal{A}(X), \mathfrak{A}(X))$ is a freely graded
quasiaffine algebra.
\end{prop}

\proof
We have to show that there is a natural pair 
$(A',I') \in \mathfrak{A}(X)$ with the properties of
Definition~\ref{freegradalgdef}.

Choose homogeneous $f_{1}, \ldots, f_{r} \in \mathcal{A}(X)$ such that
the sets $X \setminus Z(f_{i})$ form an affine cover of $X$.
Let $q \colon \rq{X} \to X$ be the canonical map.
By Proposition~\ref{geominterp}~(ii), 
each $\rq{X}_{f_{i}}$ equals $q^{-1}(X \setminus Z(f_{i}))$
and thus is affine.
Consequently the algebras
$$ 
\mathcal{A}(X)_{f_{i}} 
= 
\mathcal{O}(\rq{X})_{f_{i}} 
=
\mathcal{O}(\rq{X}_{f_{i}})
$$
are of finite type. Thus we find a
subalgebra $A' \subset \mathcal{A}(X)$ of finite type satisfying 
$A'_{f_{i}} = \mathcal{A}(X)_{f_{i}}$ for every $i$. Then 
$\b{X} := \Spec(A')$ is an affine closure of $\rq{X}$, and the
vanishing ideal $I' \subset A'$ of $\b{X} \setminus \rq{X}$ is the
radical of the ideal generated by $f_{1}, \ldots, f_{r}$. It follows
that $(A',I')$ is a natural pair on $\rq{X}$.

We verify the condition on the degrees.  Given $x \in \rq{X}$, choose an
$f_{i}$ with $q(x) \in U := X \setminus Z(f_{i})$.  By
Lemma~\ref{gradproject}~(iii), there is a small neighbourhood 
$U_{h} \subset U$ 
of $x$ defined by some $h \in \mathcal{O}(U)$ such that every 
$[L] \in \Pic(X)$ admits an invertible section in 
$\mathcal{A}_{[L]}(U_{h})$.

Now, $U_{h}$ equals $X \setminus Z(hf_{i}^{n})$ for some large
positive integer $n$.
Since finitely many of such $U_{h}$ cover $X$, we obtain the desired
Property~\ref{freegradalgdef} with finitely many of the homogeneous 
sections $hf_{i}^{n} \in I'$.
\endproof

\begin{defi}\label{homcoorddef}
We call $(\mathcal{A}(X), \mathfrak{A}(X))$ 
the {\em homogeneous coordinate ring\/} of $X$.
\end{defi}

We show now that homogeneous coordinate rings are unique up to
isomorphism. This amounts to comparing Picard graded algebras 
arising from different groups of line bundles on $X$. 
As we shall need it later, we do this in a slightly more general
setting: 

\begin{lemma}\label{differentcomps}
Let $\Lambda$ and $\Gamma$ be groups of line bundles on $X$ 
with associated graded $\mathcal{O}_{X}$-algebras $\mathcal{R}$ 
and $\mathcal{S}$.
Suppose that the image of $\Lambda \to \Pic(X)$ contains the image of
$\Gamma \to \Pic(X)$, and let $\varrho$ be a
full shifting family for $\mathcal{R}$. 
\begin{enumerate}
\item There exist a graded homomorphism $\gamma \colon \mathcal{S}
  \to \mathcal{R}$ with accompanying homomorphism $\t{\gamma} \colon
  \Gamma \to \Lambda$ and a full shifting family $\sigma$ for $\Gamma$ 
   such that for every $K \in \Gamma$ we have $K \cong \t{\gamma}(K)$,
   and, given an $F$ from the kernel of $\Gamma \to \Pic(X)$, there is
   a commutative diagram of $\mathcal{O}_{X}$-module isomorphisms
$$   
\xymatrix{
{\mathcal{S}_{K}} \ar[rr]^{\gamma_{K}} \ar[d]_{\sigma_{F}}
& &
{\mathcal{R}_{\t{\gamma}(K)}} \ar[d]^{\varrho_{\t{\gamma}(F)}} \\
{\mathcal{S}_{K+F}} \ar[rr]_{\gamma_{K+F}} 
& &
{\mathcal{R}_{\t{\gamma}(K)+\t{\gamma}(F)}}
}
$$
\item Given data as in (i), let $\mathcal{B}$ and $\mathcal{A}$ 
  denote the Picard graded algebras associated to $\sigma$ and
  $\varrho$. Then one has a commutative diagram
$$
\xymatrix{
{\mathcal{S}} \ar[r]^{\gamma} \ar[d] 
&
{\mathcal{R}} \ar[d] \\
{\mathcal{B}} \ar[r]_{\b{\gamma}} 
&
{\mathcal{A}}
}
$$
  of graded $\mathcal{O}_{X}$-algebra homomorphisms. The lower row 
  is an isomorphism if $\Gamma$ and $\Lambda$ have the same image in
  $\Pic(X)$.
\end{enumerate}
\end{lemma}

\proof 
Let $\Gamma \subset \Lambda(\mathfrak{V})$ and 
$\Lambda \subset \Lambda(\mathfrak{U})$. 
Then $\Lambda$ and $\Gamma$ embed canonically into 
$\Lambda(\mathfrak{W})$, where $\mathfrak{W}$
denotes any common refinement of the open covers 
$\mathfrak{U}$ and $\mathfrak{V}$.
Hence we may assume that $\Lambda$ and $\Gamma$ 
arise from the same trivializing cover. 

Let $K_{1}, \ldots, K_{m}$ be a basis of $\Gamma$ and choose 
$E_{1}, \ldots, E_{m} \in \Lambda$ in such a way that the isomorphism 
class of 
$E_{i}$ equals the class of $K_{i}$ in 
$\Pic(X)$. Furthermore let $\t{\gamma} \colon \Gamma \to \Lambda$ 
be the homomorphism sending $K_{i}$ to $E_{i}$. 
For each $i = 1, \ldots, m$, fix a bundle isomorphism 
$\beta_{K_{i}} \colon K_{i} \to E_{i}$.

By multiplying the local data of the these homomorphisms, 
we obtain as in the proof of Lemma~\ref{shiftfamprops}
a bundle isomorphism $\beta_{K} \colon K \to \t{\gamma}(K)$
for every $K \in \Gamma$. Shifting sections via these 
$\beta_{K}$ defines $\mathcal{O}_{X}$-module isomorphisms
$\gamma_{K} \colon \mathcal{S}_{K} \to \mathcal{R}_{\t{\gamma}(K)}$.
By construction, the $\gamma_{K}$ fit together to a graded homomorphism
$\gamma \colon \mathcal{S} \to \mathcal{R}$ of
$\mathcal{O}_{X}$-algebras. 

Now it is clear how to define the full shifting family $\sigma$: 
Take an $F$ from the kernel of $\Gamma \to \Pic(X)$. 
Define $\sigma_{F} \colon \mathcal{S} \to \mathcal{S}$ by
prescribing on the homogeneous components the (unique) isomorphisms
$\mathcal{S}_{K} \to \mathcal{S}_{K+F}$ that make the above diagrams
commutative. It is then straightforward to verify the properties of a
shifting family for the maps $\sigma_{F}$. This settles
assertion~(i).

We prove~(ii). By the commutative diagram of~(i), the ideal associated
to $\sigma$ is mapped into the ideal associated to $\varrho$. 
Hence, we obtain the desired homomorphism 
$\b{\gamma} \colon \mathcal{B} \to \mathcal{A}$
of Picard graded algebras.

Now, assume that the images of $\Gamma$ and $\Lambda$ in $\Pic(X)$
coincide. Since every $\gamma_{K} \colon \mathcal{S}_{K} \to
\mathcal{R}_{\t{\gamma}(K)}$ is an isomorphism, we can use
Lemma~\ref{gradproject}~(ii) to see that $\b{\gamma}$ is an isomorphism in
every degree. By assumption the accompanying homomorphism of $\b{\gamma}$ is
bijective, whence the assertion follows.  \endproof

The uniqueness of homogeneous coordinate rings is a direct consequence 
of the Lemmas~\ref{shiftfamunique} and~\ref{differentcomps}:

\begin{prop}\label{uniquehomcoord}
Different choices of the group of line bundles and the full shifting
family define isomorphic freely graded quasiaffine algebras as
homogeneous coordinate rings for~$X$.
\end{prop}

\proof
Let $\Lambda$ and $\Gamma$ be two groups of line bundles mapping onto
$\Pic(X)$ and let $\mathcal{A}$ and $\mathcal{B}$ denote the Picard
graded algebras associated to choices of full shifting families for
the corresponding $\Lambda$ and $\Gamma$-graded
$\mathcal{O}_{X}$-algebras. From Lemmas~\ref{shiftfamunique}
and~\ref{differentcomps} we infer the existence of a graded
$\mathcal{O}_{X}$-algebra isomorphism 
$\mu \colon \mathcal{B} \to \mathcal{A}$. 
In particular, we have $\mathcal{B}(X) \cong \mathcal{A}(X)$. 

We show that $\mu$ defines an isomorphism of quasiaffine
algebras. 
Let $(B',J') \in \mathfrak{B}(X)$ as in~\ref{freegradalgdef}. 
Then Lemma~\ref{quasiaff2closingsubalg} ensures that 
$(B',J')$ is a natural pair on $\Spec(\mathcal{B})$. 
We have to show that $(\mu(B'), \mu(I'))$ is a natural pair on
$\Spec(\mathcal{A})$. 
Since $\mu$ is an $\mathcal{O}_{X}$-module 
isomorphism in every degree, we have $Z(\mu(g)) = Z(g)$
for any homogeneous $g \in J'$.
Thus Proposition~\ref{geominterp}~(ii)
tells us that $(\mu(B'), \mu(I'))$ is a natural pair.
\endproof


\section{Functoriality of the homogeneous coordinate ring}%
\label{section5}

In this section, we present the first main result. It
says that homogeneous coordinates are a fully faithful 
contravariant functor, see Theorem~\ref{fullyfaithful}.
But first we have to define the homogeneous coordinate ring
functor on morphisms.
The basic tool for this definition are Picard graded 
pullbacks, see~\ref{picpulldef} and~\ref{picpullex}.

\goodbreak

As in the preceding section, we assume that $\KK^{*}$ is 
of infinite rank over $\ZZ$. 
Moreover, in this section we assume all varieties to be
divisorial and to have only constant
invertible global functions.  
Finally, we require that any variety has a 
finitely generated Picard group,
and, if $\KK$ is of characteristic $p > 0$,
this Picard group has no $p$-torsion.

For a variety $X$ fix a group $\Lambda$ of line bundles mapping onto $\Pic(X)$
and denote the associated $\Lambda$-graded $\mathcal{O}_{X}$-algebra by
$\mathcal{R}$. 
Moreover, we fix a full shifting family $\varrho$ for
$\mathcal{R}$ and denote the resulting Picard graded algebra by
$\mathcal{A}$. 
For a further variety $Y$ we denote the corresponding data by
$\Gamma$, $\mathcal{S}$, $\sigma$, and $\mathcal{B}$.  
Let $\varphi \colon X \to Y$ 
be a morphism of the varieties $X$ and $Y$.

\begin{defi}\label{picpulldef}
By a {\em Picard graded pullback for $\varphi \colon X \to Y$} 
we mean a graded homomorphism 
$\mathcal{B} \to \varphi_{*}\mathcal{A}$
of $\mathcal{O}_{Y}$-algebras having the pullback map 
$\varphi^{*} \colon \Pic(Y) \to \Pic(X)$ as its accompanying 
homomorphism.
\end{defi}

Note that the property of being an $\mathcal{O}_{Y}$-algebra
homomorphism means in particular, that in degree zero any Picard
graded pullback is the usual pullback of functions. 
As a consequence, we remark:

\begin{lemma}\label{pullzero}
Let $\mu \colon \mathcal{B} \to \varphi_{*}\mathcal{A}$ be a 
Picard graded pullback for $\varphi \colon X \to Y$, and let $g \in
\mathcal{B}(Y)$ be homogeneous. 
Then the zero set $Z(\mu(g)) \subset X$ is the inverse 
image $\varphi^{-1}(Z(g))$ of the zero set $Z(g) \subset Y$.
\end{lemma}

\proof 
It suffices to prove the statement locally, 
over small open $V \subset Y$.
But on such $V$, we may shift $g$ by multiplication with 
invertible elements into degree zero. 
This does not affect zero sets, 
whence the assertion follows.
\endproof

The basic step in the definition of the homogeneous coordinate
ring functor on morphisms is to show existence of Picard graded
pullbacks and to provide a certain uniqueness property:

\begin{prop}\label{picpullex}
There exist Picard graded pullbacks for $\varphi \colon X \to Y$. 
Moreover, any two Picard graded pullbacks 
$\mu, \nu \colon \mathcal{B} \to \varphi_{*}\mathcal{A}$
for $\varphi$ differ by a character 
$c \colon \Pic(Y) \to \KK^{*}$ 
in the sense that $\nu_{P} = c(P) \mu_{P}$
holds for all $P \in \Pic(Y)$.
\end{prop}

The proof of this statement is based on two lemmas. The first one
is an extension property for shifting families:

\begin{lemma}\label{shiftext}
Let $\Pi$ be any group of line bundles on $X$, and let $\Pi_{0}
\subset \Pi_{1}$ be two subgroups of the kernel of $\Pi \to \Pic(X)$. 
Then every $\Pi_{0}$-shifting family $\tau^{0}$ for the $\Pi$-graded
$\mathcal{O}_{X}$-algebra $\mathcal{T}$ associated to $\Pi$
extends to a $\Pi_{1}$-shifting family $\tau^{1}$ for $\mathcal{T}$ 
in the sense that $\tau^{1}_{E} = \tau^{0}_{E}$ holds for all
$E \in \Pi_{0}$. 
\end{lemma}

\proof
Let $\vartheta$ be any $\Pi_{1}$-shifting family for $\mathcal{T}$.
Then $\vartheta$ restricts to a $\Pi_{0}$-shifting family. By
Lemma~\ref{shiftfamprops}, there is a character 
$c \colon \Pi_{0} \to \mathcal{O}^{*}(X)$ 
with $\tau^{0}_{E} = c(E) \vartheta_{E}$ for all $E \in \Pi_{0}$. 
As we assumed $\mathcal{O}^{*} (X) = \KK^{*}$, Lemma~\ref{charext} tells us 
that $c$ extends to $\Pi_{1}$. Thus, setting
$\tau^{1}_{E} := c(E) \vartheta_{E}$ for $E \in \Pi_{1}$ gives
the desired extension.
\endproof

The second lemma provides a pullback construction for shifting
families. By pulling back cocycles, we obtain the (again free) 
pullback group $\varphi^{*}\Gamma$. We denote the associated
$\varphi^{*}\Gamma$-graded $\mathcal{O}_{X}$-algebra by
$\varphi^{*}\mathcal{S}$. Indeed $\varphi^{*}\mathcal{S}$ 
is canonically isomorphic to the ringed inverse image of 
$\mathcal{S}$. 
Observe that we have a canonical sheaf homomorphism 
$\mathcal{S} \to \varphi_{*}\varphi^{*} \mathcal{S}$.

\goodbreak

\begin{lemma}\label{pullshift}
Let $\Gamma_{0} \subset \Gamma$ a subgroup, and let
$\sigma$ be a $\Gamma_{0}$-shifting family for $\mathcal{S}$.
\begin{enumerate}
\item The $\mathcal{O}_{X}$-module homomorphisms
  $\varphi^{*}\sigma_{F}$ define a
  $\varphi^{*}\Gamma_{0}$-shifting family
  $\varphi^{*}\sigma$ for $\varphi^{*}\mathcal{S}$.
\item The ideal $\mathcal{J}^{*}$ associated to $\varphi^{*}\sigma$
  equals the pullback $\varphi^{*}\mathcal{J}$ of the ideal
  $\mathcal{J}$ associated to $\sigma$.
\end{enumerate} 
\end{lemma}

\proof
For~(i), note that the isomorphisms $\sigma_{F} \colon \mathcal{S}_{K}
\to \mathcal{S}_{K+F}$ can be written as $g \mapsto \beta_{K,F}(g)$
with unique line bundle isomorphisms $\beta_{K,F} \colon K \to K + F$.
The family $\varphi^{*} \sigma_{F}$ corresponds to the collection
$\varphi^{*} \beta_{K,F} \colon \varphi^{*} K \to \varphi^{*}K +
\varphi^{*}F$. The properties of a shifting family become clear
by writing the $\varphi^{*} \beta_{K,F}$ in terms of local data
as in~\ref{localdata}.

To prove~(ii), we just compare the stalks of the two sheaves in
question. By Property~\ref{shiftfamprops}~(iii), we obtain for any $x
\in X$:
\begin{eqnarray*}
\mathcal{J}^{*}_{x}
& = & 
\bangle{1_{x} - \varphi^{*} \sigma_{F}(1_{x}); 
  \; F \in \Gamma_{0}} \\
& = & 
\bangle{\varphi^{*}(1_{\varphi(x)}) - 
 \varphi^{*}(\sigma_{F}(1_{\varphi(x)})); 
  \; F \in \Gamma_{0}} \\
& = & 
(\varphi^{*}\mathcal{J})_{x}. \qquad \qed
\end{eqnarray*}

\proof[Proof of Proposition~\ref{picpullex}]
We establish the existence of Picard graded pullbacks:
As usual, let $\mathcal{I}$ and $\mathcal{J}$ 
denote the respective ideals associated  to the shifting families
$\varrho$ for $\mathcal{R}$ and $\sigma$ for $\mathcal{S}$. 
Thus the corresponding Picard graded algebras are 
$\mathcal{A} = \mathcal{R}/\mathcal{I}$ and 
$\mathcal{B} = \mathcal{S}/\mathcal{J}$.

By Lemma~\ref{pullshift}, we
have the $\varphi^{*}\Gamma_{0}$-shifting family $\varphi^{*}\sigma$ 
for $\varphi^{*}\mathcal{S}$.  
Lemma~\ref{shiftext} enables us to choose a full shifting family
$\varphi^{\sharp}\sigma$ extending $\varphi^{*}\sigma$.
We denote by $\varphi^{\sharp}\mathcal{J}$ the ideal associated to 
$\varphi^{\sharp}\sigma$, and write 
$\varphi^{\sharp}\mathcal{B} := 
\varphi^{*}\mathcal{S}/\varphi^{\sharp}\mathcal{J}$
for the quotient.
In this notation, we have a commutative diagram of graded
$\mathcal{O}_{Y}$-algebra homomorphisms such that the unlabelled
arrows are isomorphisms in each degree:
$$
\xymatrix{
{\varphi_{*}\mathcal{R}} \ar[d]
& & 
{\varphi_{*}\varphi^{*}\mathcal{S}} \ar[dl] \ar[dr]  \ar[ll]
& &
{\mathcal{S}} \ar[ll]_{\varphi^{*}} \ar[d] \\
{\varphi_{*}{\mathcal{A}}} 
&
{\varphi_{*}\varphi^{\sharp}\mathcal{B}} \ar[l]
&
&
{\varphi_{*}\varphi^{*}\mathcal{B}} \ar[ll]
&
{\mathcal{B}} \ar[l]^{\varphi^{*}}
}
$$

Indeed, the right square is standard.  To obtain the
middle triangle, we only have to show that $\varphi^{\sharp}\mathcal{J}$
contains the kernel of $\varphi^{*}\mathcal{S} \to
\varphi^{*}\mathcal{B}$. But this follows from exactness of $\varphi^{*}$ and
Lemma~\ref{pullshift}~(ii). Existence of the left square follows from
combining Lemmas~\ref{shiftfamunique} and~\ref{differentcomps}.
Now the desired Picard graded pullback of $\varphi \colon X \to Y$ is
the composition of the lower horizontal arrows.

We turn to the uniqueness statement.
Let $P \in \Pic(Y)$. Since $\mathcal{B}_{P}$ locally admits
invertible sections, we can cover $Y$ by open $V \subset Y$ such that
there exist invertible sections $h \in \mathcal{B}_{P} (V)$. We define 
$$c(P,V) := \nu(h)/\mu (h) \in \mathcal{A}_{0}^{*} (\varphi^{-1}(V)).$$ 

This does not depend on the choice of $h$: For a further invertible 
$ g\in \mathcal{B}_{P} (V)$, the section $g/h$ is of degree zero.
But in degree zero any Picard graded pullback is the usual pullback 
of functions. 
Thus we have $\mu(h/g) = \nu(h/g)$.
Consequently, $\nu(g)/\mu(g)$ equals $\nu(h)/\mu(h)$.

\goodbreak

Similarly we see that for two open $V, V' \subset Y$ as above, the
corresponding sections $c(P,V)$ and $c(P,V')$ coincide on the
intersection $\varphi^{-1}(V) \cap \varphi^{-1}(V')$.
Thus, by gluing, we obtain a global
section $c(P) \in \mathcal{A}_{0}^{*} (X) = \mathcal{O}^{*}(X)$. 
Then it is immediate to check, that $P \mapsto c(P)$ has the
desired properties.
\endproof

With the help of Picard graded pullbacks we can now 
make the homogeneous coordinate ring into a functor.
We fix for any morphism $\varphi \colon X \to Y$ 
a Picard graded pullback
$\mu_{\varphi} \colon \mathcal{B} \to \varphi_{*}\mathcal{A}$,
and denote the induced homomorphism on global sections again by 
$\mu_{\varphi} \colon\mathcal{B}(Y) \to \mathcal{A}(X)$.

For a graded homomorphism $\nu$ of freely graded quasiaffine
algebras, we denote by $[\nu]$ its equivalence class in the sense of
Definition~\ref{pointedmorphdef}.

\begin{prop}\label{var2alghom}
The assignments $X \mapsto (\mathcal{A}(X), \mathfrak{A}(X))$
and $\varphi \mapsto [\mu_{\varphi}]$ define a contravariant
functor into the category of freely graded quasiaffine algebras.
\end{prop}

\proof 
By Proposition~\ref{coringisqualg},
the homogeneous coordinate rings 
$(\mathcal{A}(X),\mathfrak{A}(X))$ 
and $(\mathcal{B}(Y),\mathfrak{B}(Y))$ 
of $X$ and $Y$ are
in fact freely graded quasiaffine algebras. 
The first task is to show that the homomorphism
$\mu_{\varphi} \colon \mathcal{B}(Y) \to \mathcal{A}(X)$ 
associated to a morphism $\varphi \colon X \to Y$
is a graded homomorphism of freely graded quasiaffine 
algebras.

As a Picard graded pullback, $\mu_{\varphi}$ is graded and 
has as accompanying homomorphism the pullback map 
$\Pic(Y) \to \Pic(X)$. 
Thus we are left with checking the conditions of
Definition~\ref{quasiaffalgdef}~(ii) for $\mu_{\varphi}$. 
This is done geometrically in terms of the constructions of
Proposition~\ref{geominterp}:
$$ \rq{X} := \Spec(\mathcal{A}), \qquad 
\rq{Y} := \Spec(\mathcal{B}), \qquad  
 q_{X} \colon \rq{X} \to X, \qquad q_{Y} \colon \rq{Y} \to Y. $$

Let $(B',J') \in \mathfrak{B}(Y)$ be a closing subalgebra as in
Definition~\ref{freegradalgdef}. Then Lemma~\ref{naturalpairs}
provides a closing subalgebra $(A',I') \in \mathfrak{A}(X)$ 
such that $\mu_{\varphi}(B') \subset A'$ holds.
We have to verify the condition on the ideals $I'$ and $J'$
required in~\ref{quasiaffalgdef}~(ii).  
For this, consider the affine closures
of $\rq{X}$ and $\rq{Y}$:
$$\b{X} := \Spec(A'), \qquad \b{Y} := \Spec(B'). $$

Then the restricted homomorphism $\mu_{\varphi} \colon B' \to A'$ 
defines a morphism $\b{\varphi} \colon \b{X} \to \b{Y}$. 
Recall from Section~\ref{section1} that $I'$ and $J'$ are 
the vanishing ideals of the complements $\b{X} \setminus \rq{X}$
and $\b{Y} \setminus \rq{Y}$.
Thus we have to show that $\b{\varphi}$ maps $\rq{X}$ to $\rq{Y}$. 
For this, let $g_{1}, \ldots, g_{s} \in J'$ be
homogeneous sections as in~\ref{freegradalgdef}.
Using Lemma~\ref{pullzero}, we obtain: 

\goodbreak

\begin{eqnarray*}
\rq{X} & = &
\bigcup_{j=1}^{r} q_{X}^{-1}(\varphi^{-1}(Y \setminus Z(g_{j}))) \\
& = &
\bigcup_{j=1}^{r} q_{X}^{-1} (X \setminus Z(\mu_{\varphi}(g_{j}))) \\
& \subset & 
\bigcup_{j=1}^{r} \b{X}_{\mu_{\varphi}(g_{j})} \\
& = & 
\b{\varphi}^{-1}(\rq{Y}).
\end{eqnarray*}

\goodbreak

Finally, we check that $\varphi \mapsto [\mu_{\varphi}]$ is functorial.
Note that by Proposition~\ref{picpullex}, the class $[\mu_{\varphi}]$
does not depend on the choice of the Picard graded pullback 
$\mu_{\varphi}$ of a given morphism. 

From this we conclude that the identity morphism of a variety is
mapped to the identity 
morphism of its homogeneous coordinate ring. 
Moreover, as the composition of two Picard graded
pullbacks is a Picard graded pullback for the composition of the respective
morphisms, the above assignment commutes with composition.
\endproof

In the sequel we shall speak of the homogeneous coordinate ring functor.
We present the first main result of this article. 
It tells us that the morphisms of two varieties are in one-to-one 
correspondence with the morphisms of their coordinate rings:

\begin{thm}\label{fullyfaithful}
The homogeneous coordinate ring functor $X \mapsto (\mathcal{A}(X),
\mathfrak{A}(X))$ and $\varphi \mapsto [\mu_{\varphi}]$ is fully faithful.
\end{thm}

\proof Let $X$, $Y$ be varieties with associated Picard graded
algebras $\mathcal{A}$ and $\mathcal{B}$. 
We denote the respective homogeneous
coordinate rings of $X$ and $Y$ for short by $(A,\mathfrak{A})$ and
$(B,\mathfrak{B})$. We construct an inverse to
$$\Mor(X,Y) \to \Mor((B,\mathfrak{B}),(A,\mathfrak{A})),
\qquad \varphi \mapsto [\mu_{\varphi}].$$

So, start with any graded homomorphism $\mu \colon (B,\mathfrak{B}) \to
(A,\mathfrak{A})$ of quasiaffine algebras. Then Lemma~\ref{naturalpairs}
provides closing subalgebras $(A',I') \in \mathfrak{A}$ and $(B',J') \in
\mathfrak{B}$ such that $(B',J')$ is as in Definition~\ref{freegradalgdef}
and we have $\mu(B') \subset A'$.

Consider the affine closures $\b{X} := \Spec(A')$ and
$\b{Y} := \Spec(B')$ of $\rq{X} := \Spec(\mathcal{A})$ and $\rq{Y} :=
\Spec(\mathcal{B})$. Then $\mu$ gives rise to a morphism $\b{\varphi}
\colon \b{X} \to \b{Y}$, and restricting this morphism to $\rq{X}$
yields a commutative diagram 
$$ \xymatrix{%
{\rq{X}}  
\ar[r]^{\rq{\varphi}} 
\ar[d]_{q_{X}} & 
{\rq{Y}} 
\ar[d]^{q_{Y}}  \\
X
\ar[r]^{\varphi} & 
Y } $$
where $q_{X}$ and $q_{Y}$ denote the canonical maps, and 
the morphism $\varphi \colon X \to Y$ has as its pullbacks on the
level of functions the maps obtained by restricting 
the localizations $\mu_{g} \colon B_{g} \to A_{\mu(g)}$ to degree zero
over the affine sets $Y_{g}$ for homogeneous $g \in J'$. 

Observe that applying the above procedure 
to a further graded homomorphism 
$\nu \colon (B,\mathfrak{B}) \to (A,\mathfrak{A})$ yields the same
induced morphism $X \to Y$ if and only if the homomorphisms $\mu$ and
$\nu$ are equivalent; the ``only if'' part follows from the uniqueness
statement of Proposition~\ref{picpullex} and the fact that $\mu$ and
$\nu$ define Picard graded pullbacks via localizing. Thus
$[\mu] \mapsto \varphi$ defines an injection
$$\Mor((B,\mathfrak{B}),(A,\mathfrak{A})) \to \Mor(X,Y).$$

We check that this map is inverse to the one defined by the
homogeneous coordinate ring functor. Start with a morphism
$\varphi \colon X \to Y$, and let 
$[\mu_{\varphi}] \colon B \to A$ be as 
before Proposition~\ref{var2alghom}. 
Write shortly $\mu := \mu_{\varphi}$.
Consider a homogeneous $g \in B$
such that $V := Y \setminus Z(g)$ is affine and let $U :=
\varphi^{-1}(V)$. 
Using Lemma~\ref{pullzero}, we obtain a commutative diagram 
$$\xymatrix{
{A_{(\mu(g))}} \ar@{=}[d] &
{B_{(g)}}  \ar[l]_{\mu_{(g)}} \ar@{=}[d] \\
{\mathcal{O}_{X}(U)}  & 
{\mathcal{O}_{Y}(V)} \ar[l]^{\varphi^{*}}   
}
$$
where the above horizontal map is the map on degree zero induced by
the localized map $\mu_{g} \colon B_{g} \to A_{\mu(g)}$. 
Since $Y$ is covered by open affine sets of the form 
$V = Y \setminus Z(g)$, 
we see that the morphism $X \to Y$
associated to $\mu = \mu_{\varphi}$ 
is again $\varphi$. \endproof 

So far, our homogeneous coordinate ring functor depends on the choice 
of the homogeneous coordinate ring for a given variety.
By passing to isomorphism classes, the whole
construction can even be made canonical:

\begin{rem}\label{canonical}
If one takes as target category the category of isomorphism classes of
freely graded quasiaffine algebras, then the homogeneous coordinate
functor $X \to (\mathcal{A}(X), \mathfrak{A}(X))$ and $\varphi \mapsto
[\mu_{\varphi}]$ becomes unique. 
\end{rem}

\section{A first dictionary}\label{section6}

We present a little dictionary between geometric
properties of a variety and algebraic 
properties of its homogeneous coordinate ring. 
We consider separatedness, normality and smoothness.
Moreover, we treat quasicoherent sheaves, 
and we describe affine morphisms and closed embeddings.

The setup is the same as in Sections~\ref{section4}
and~\ref{section5}:
The multiplicative group $\KK^{*}$ 
of the ground field $\KK$ is supposed to 
be of infinite rank over $\ZZ$. 
Moreover, $X$ is a divisorial variety with 
$\mathcal{O}^{*}(X) = \KK^{*}$ and its Picard group
is finitely generated and has no $p$-torsion if $\KK$ is of characteristic
$p>0$.

Denote by $(A,\mathfrak{A}) := (\mathcal{A}(X),\mathfrak{A}(X))$ 
the homogeneous coordinate ring of $X$. 
Recall that $A$ is the algebra of global sections of a suitable 
Picard graded $\mathcal{O}_{X}$-algebra~$\mathcal{A}$.
In the subsequent proofs, we shall often use
the geometric interpretation provided by
Propositions~\ref{geominterp} and~\ref{equivquasiaffequiv}:

\begin{lemma}\label{geomquot}
Consider $\rq{X} := \Spec(\mathcal{A})$, the canonical map 
$q \colon \rq{X} \to X$ and the diagonalizable group
$H := \Spec(\KK[\Pic(X)])$.
\begin{enumerate}
\item There is a unique free action of $H$ on
  $\rq{X}$ such that each $\mathcal{A}_{[L]}(U)$ consists precisely
  of the $\chi^{[L]}$-homogeneous functions of $q^{-1}(U)$.
\item The canonical map $q \colon \rq{X} \to X$ is a geometric
  quotient for the above $H$-action on $X$. 
\end{enumerate}
\end{lemma}

\proof The first statement follows from Propositions~\ref{geominterp}
and~\ref{equivquasiaffequiv}. The second statement is due to the facts
that $\mathcal{O}_{X} = q_{*}(\mathcal{A}_{0})$ is the sheaf of
invariants and the action of $H$ is free. \endproof

We begin with the dictionary. It is quite easy to characterize
separatedness in terms of the homogeneous coordinate ring:

\begin{prop}
The variety $X$ is separated if and only if there exists a 
graded closing subalgebra $(A',I') \in \mathfrak{A}$ and 
homogeneous $f_{1}, \ldots, f_{r} \in I'$ 
as in~\ref{freegradalgdef} such that each of the maps
$A_{(f_{i})} \otimes A_{(f_{j})} \to A_{(f_{i}f_{j})}$
is surjective. 
\end{prop}

\proof 
First recall that the sets 
$X_{i} := X \setminus Z(f_{i})$ 
form an affine cover of $X$. 
The above condition means just that 
the canonical maps from
$\mathcal{O}(X_{i}) \otimes \mathcal{O}(X_{j})$
to $\mathcal{O}(X_{i} \cap X_{j})$ 
is surjective. 
This is the usual separatedness criterion~\cite[Prop.~3.3.5]{Ke}. 
\endproof

Next we show how normality of the variety $X$ is
reflected in its homogeneous coordinate ring 
(for us, a normal variety is in particular
irreducible):

\begin{prop}\label{normal2normal}
The variety $X$ is normal if and only if $A$ 
is a normal ring.
\end{prop}

\proof
We work in terms of the geometric data 
$q \colon \rq{X} \to X$ and $H$ discussed in Lemma~\ref{geomquot}.
First suppose that $A = \mathcal{A}(X)$ is a normal ring. 
Then the quasiaffine variety $\rq{X}$ is normal. 
It is a basic property of geometric quotients
that the variety $X$ inherits normality from $\rq{X}$,
see e.g.~\cite[p.~39]{Do}.

Suppose conversely that $X$ is normal. 
Luna's Slice Theorem
tells us that $q \colon \rq{X} \to X$ 
is an $H$-principal bundle in the \'etale topology,
see~\cite{Lu}, and~\cite[Prop.~8.1]{BaRi}.
Thus, up to \'etale maps, $\rq{X}$ looks locally 
like $X \times H$.
Since normality of local rings is stable under \'etale 
maps~\cite[Prop.~I.3.17]{Mi}, 
we can conclude that all local rings of $\rq{X}$ are 
normal. 

It remains to show that $\rq{X}$ is connected.
Assume the contrary. Then there is a connected
component $\rq{X}_{1} \subset \rq{X}$ with $q(\rq{X}_{1}) = X$. 
Let $H_{1} \subset H$ be the stabilizer of $\rq{X}_{1}$, that
means that $H_{1}$ is the maximal subgroup of $H$
with $H_{1} \mal \rq{X}_{1} = \rq{X}_{1}$.
Note that we have $t \in H_{1}$ if $t \mal x \in \rq{X}_{1}$ holds 
for at least one point $x \in \rq{X}_{1}$. In particular, $H_{1}$ is a 
proper subgroup of $H$.

We claim that restricting the canonical map $q \colon \rq{X} \to X$ 
to $\rq{X}_{1}$ yields a geometric quotient for the action of $H_{1}$
on $\rq{X}_{1}$. Indeed, $H_{1}$ acts freely on $\rq{X}_{1}$. Hence
we have a geometric quotient $\rq{X}_{1} \to \rq{X}_{1}/H_{1}$ and a 
commutative diagram
$$ \xymatrix{ 
{\rq{X}_{1}} \ar[r]^{\subset} \ar[d]_{/H_{1}} 
& 
{\rq{X}} \ar[d]^{/H}_{q} \\
{\rq{X}_{1}/H_{1}} \ar[r]
&
X }
$$

The map $\rq{X}_{1}/H_{1} \to X$ is bijective, because the 
intersection of a $q$-fibre with~$\rq{X}_{1}$ always is 
a single $H_{1}$-orbit. 
Since $X$ is normal, we may apply Zariski's Main Theorem 
to conclude that $\rq{X}_{1}/H_{1} \to X$ is even
an isomorphism. This verifies our claim. 

Since $H_{1}$ is a proper subgroup of $H$,
we find a nontrivial class $[L] \in \Pic(X)$ such that the
corresponding character $\chi^{[L]}$ of $H$ is 
trivial on $H_{1}$.
We construct a defining cocycle for the class $[L]$: Cover $X$ by
small open sets $U_{i}$ admitting invertible sections 
$g_{i} \in \mathcal{A}_{[L]}(U_{i})$. Then the cocycle $g_{i}/g_{j}$ defines
a bundle belonging to the class $[L]$.

On the other hand, the $g_{i}$ are $\chi^{[L]}$-homogeneous 
functions on $q^{-1}(U_{i})$. So they restrict to 
$H_{1}$-invariant functions on $q^{-1}(U_{i}) \cap \rq{X}_{1}$. 
As seen before, $X$ is the quotient 
of $\rq{X}_{1}$ by the action of $H_{1}$. 
Thus we conclude that the $g_{i}/g_{j}$ form in fact a 
coboundary on $X$. 
Consequently, the class $[L]$ must be trivial. 
This contradicts the choice of $[L]$.
\endproof

Thus we see that if $X$ is normal, then
$A$ is the ring of global functions of a 
normal variety. 
That means that $A$ belongs to a intently studied 
class of rings:

\begin{coro}
Let $X$ be normal. Then $A$ is a Krull ring. \endproof
\end{coro}

As we did in Proposition~\ref{normal2normal} for normality,
we can characterize smoothness in terms of the homogeneous coordinate
ring:

\begin{prop}\label{smooth}
$X$ is smooth if and only if there is  
a closing subalgebra $(A',I') \in \mathfrak{A}$ 
such that all localizations
$A_{\mathfrak{m}}$ are regular, where $\mathfrak{m}$ 
runs through the maximal ideals with 
$I' \not \subset \mathfrak{m}$.
\end{prop}

\proof 
Let $\rq{X} := \Spec(\mathcal{A})$, 
and consider the affine closure $\b{X} := \Spec(A')$
defined by any closing subalgebra $(A',I')$
of $A$.
Recall from Lemmas~\ref{naturalpairs} and~\ref{quasiaff2closingsubalg}
that $I'$ is the vanishing ideal
of the complement $\b{X} \setminus \rq{X}$.
So, the regularity of the local rings
$A_{\mathfrak{m}}$, where $I' \not \subset \mathfrak{m}$,
just means smoothness of $\rq{X}$.

The rest is similar to the proof of Proposition~\ref{normal2normal}:
The canonical map $q \colon \rq{X} \to X$
is an \'etale $H$-principal bundle for a diagonalizable 
group $H$.
Thus, up to \'etale maps, $\rq{X}$ looks locally 
like $X \times H$.
The assertion then follows from the fact 
that regularity of local rings is 
stable under \'etale maps, see~\cite[Prop.~I.3.17]{Mi}.
\endproof

We give a description of quasicoherent sheaves.
Consider a graded $A$-module $M$. Given $f_{1}, \ldots, f_{r} \in
A$ as in~\ref{freegradalgdef}, set $\mathcal{M}_{i} :=
M_{(f_{i})}$. Then these modules glue together to a quasicoherent
$\mathcal{O}_{X}$-module $\mathcal{M}$ on $X$. As in the toric case
\cite[Section~4]{ACHaSc}, one obtains:

\begin{prop}
The assignment $M \mapsto \mathcal{M}$ defines an essentially
surjective functor from the category of graded $A$-modules to the
category of quasicoherent $\mathcal{O}_{X}$-modules. \endproof
\end{prop}

We come to properties of morphisms.
Let $Y$ be a further variety like $X$, 
and denote its homogeneous coordinate ring by $(B,\mathfrak{B})$. 
Let $\varphi \colon X \to Y$ be any morphism. 
Denote by 
$[\mu] \colon (B,\mathfrak{B}) \to (A,\mathfrak{A})$ 
the corresponding morphism of freely graded 
quasiaffine algebras. 

\begin{prop}
The morphism $\varphi \colon X \to Y$ is affine if and only if there
are graded closing subalgebras $(A',I') \in \mathfrak{A}$ and $(B',J')
\in \mathfrak{B}$ satisfying~\ref{quasiaffalgdef} such that 
$$ \sqrt{I'} = \sqrt{\bangle{\mu(J')}}. $$
Moreover, $\varphi \colon X \to Y$ is a closed embedding if and
only if it satisfies the above condition and, given $g_{1}, \ldots,
g_{s} \in B$ as in \ref{freegradalgdef}, every induced map
$B_{(g_{i})} \to A_{(\mu(g_{i}))}$ is surjective.
\end{prop}

\proof 
Let $\mathcal{B}$ be a Picard graded 
$\mathcal{O}_{Y}$-algebra with
$B = \mathcal{B}(Y)$.
Consider the affine closures $\b{X} := \Spec(A')$ 
of $\rq{X} := \Spec(\mathcal{A})$ 
and $\b{Y} := \Spec(B')$ of 
$\rq{Y} := \Spec(\mathcal{B})$. 
Then $\mu \colon B
\to A$ gives rise to a commutative diagram 
$$ \xymatrix{%
{\rq{X}}  
\ar[r]^{\rq{\varphi}} 
\ar[d]_{q_{X}} & 
{\rq{Y}} 
\ar[d]^{q_{Y}}  \\
X
\ar[r]^{\varphi} & 
Y } $$
The morphism $\varphi$ is affine if and only if $\rq{\varphi}$ is
affine. The latter is equivalent to the condition of 
$\sqrt{I'} = \sqrt{\bangle{\mu(J')}}$ of the assertion. 
The supplement on embeddings is obvious. \endproof

\section{Tame varieties}\label{section7}

In this section we shed some light on the question 
which freely graded quasiaffine algebras
occur as homogeneous coordinate rings.
As before, we assume that the multiplicative group
$\KK^{*}$ of the ground field is of infinite rank over 
$\ZZ$. We consider varieties of the following type:

\begin{defi}
A {\em tame variety\/} is a normal divisorial 
variety $X$ with $\mathcal{O}(X) = \KK$ and 
a finitely generated Picard group
$\Pic(X)$ having no $p$-torsion if $\KK$ is 
of characteristic $p>0$.
\end{defi}

The prototype of a tame variety lives in characteristic zero,
and is a smooth complete variety with finitely generated 
Picard group.
Moreover, in characteristic zero, 
every Calabi-Yau variety is tame,
and every $\QQ$-factorial rational variety $X$ with 
$\mathcal{O}(X) = \KK$ is tame.  
Finally, in characteristic zero every normal divisorial 
variety with finitely generated Picard group admits 
an open embedding into a tame variety. 

In order to figure out the coordinate rings of tame varieties,
we need some preparation. Suppose that an algebraic group $G$ acts 
on a variety $X$. Recall that a {\em $G$-linearization} of a line
bundle $E \to X$ is a fibrewise linear $G$-action on $E$ making the 
projection equivariant. 
By a {\em simple $G$-variety} we mean a $G$-variety for which any 
$G$-linearizable line bundle is trivial.

\begin{defi}
Let $\Lambda$ be a finitely generated abelian group, 
and let $(A,\mathfrak{A})$ be a freely $\Lambda$-graded 
quasiaffine algebra.
\begin{enumerate}
\item We say that $(A,\mathfrak{A})$ is {\em pointed\/} if 
  $A$ is a normal ring, $A_{0} = \KK$ holds, 
  and the set $A^{*} \subset A$ of
  invertible elements is just $\KK^{*}$.
\item We say that $(A,\mathfrak{A})$ is {\em simple\/} 
  if $\Lambda$ has no $p$-torsion if $\KK$ is 
  of characteristic $p>0$, and the quasiaffine 
  $\Spec(\KK[\Lambda])$-variety corresponding to 
  $(A,\mathfrak{A})$ is simple.
\end{enumerate}
\end{defi}

These two subclasses define full subcategories of the categories of
divisorial varieties with finitely generated Picard group and freely
graded quasiaffine algebras. The second main result of this article is
the following:

\begin{thm}\label{equivthm}
The homogeneous coordinate ring functor restricts to an
equivalence from the category of tame varieties to the category of
simple pointed algebras. 
\end{thm}

\proof
Let $X$ be a tame variety with Picard group $\Pi := \Pic(X)$, 
and denote the associated homogeneous coordinate ring
by $(A,\mathfrak{A})$. 
Then $A$ is the algebra of global sections of some 
Picard graded $\mathcal{O}_{X}$-algebra $\mathcal{A}$ on $X$. 
We shall use again the geometric data discussed in 
Lemma~\ref{geomquot}:
$$ 
\rq{X} := \Spec(\mathcal{A}),
\qquad 
q \colon \rq{X} \to X,
\qquad
H := \Spec(\KK[\Pi]). 
$$

The first task is to show that $(A,\mathfrak{A})$ 
is in fact pointed. 
From Proposition~\ref{normal2normal} we infer that 
$A$ is a normal ring.
Since we assumed $\mathcal{O}(X) = \KK$, 
and $\mathcal{O}(X)$ equals $A_{0}$, we have
$A_{0} = \KK$. 
So we have to verify $A^{*} = \KK^{*}$. 
For this, consider an arbitrary element 
$f \in A^{*}$.

Choose a direct decomposition of $\Pi$ into a free part 
$\Pi_{0}$ and the torsion part $\Pi_{{\rm t}}$. 
This corresponds to a splitting $H = H_{0} \times H_{{\rm t}}$ 
with a torus $H_{0}$ and a finite group $H_{{\rm t}}$.
As an invertible element of $\mathcal{O}(\rq{X})$, 
the function $f$ is necessarily $H_{0}$-homogeneous, 
see~e.g.~\cite[Prop.~1.1]{Ma}.
Thus, there is a degree $P \in \Pi_{0}$ such that
$$ 
f = \sum_{G \in \Pi_{t}} f_{P+G},
\qquad
f^{-1}
= \sum_{G \in \Pi_{t}} f^{-1}_{-P+G}. $$

From the identity $ff^{-1} = 1$ we infer that 
$f_{P+G}f^{-1}_{-P-G} \ne 0$ holds for at least one 
component $f_{P+G}$ of $f$.
Since $\mathcal{O}(X)=\KK$ holds, we see that 
the homogeneous section $f_{P+G} \in A$ 
is invertible. 
Thus the homogeneous component $\mathcal{A}_{P+G}$ 
is isomorphic to~$\mathcal{O}_{X}$.

On the other hand we noted in~\ref{gradproject}~(ii)
that $\mathcal{A}_{P+G}$ is isomorphic to the sheaf of 
sections of a bundle representing the class $P+G$ in 
$\Pi_{0} \oplus \Pi_{{\rm t}}$.
Thus $P+G$ is trivial, and we obtain $P = 0$.
Hence all homogeneous components of $f$ have torsion degree. 
By $\mathcal{O}(X) = \KK$ this yields that $f_{G} =0$ 
if $G \ne 0$.
Thus we have $f \in A_{0} = \KK$.

The next task is to show that $\rq{X}$ is a simple $H$-variety.
For this, let $\Pic_{H}(\rq{X})$ denote the group of equivariant
isomorphy classes of $H$-linearized line bundles on~$\rq{X}$,
compare~\cite[Sec.~2]{Kr}.
Moreover, let $\Piclin(\rq{X}) \subset \Pic(\rq{X})$ denote the
subgroup of the classes of all $H$-linearizable bundles.
We have to show that $\Piclin(\rq{X})$ is trivial. 

First, we consider the possible linearizations of the trivial bundle
$\rq{X} \times \KK$. 
Using $\mathcal{O}^{*}(\rq{X}) = \KK^{*}$, as verified before, 
one directly checks that any linearization of the trivial bundle
is given by a character $\chi$ of $H$ as follows:
\begin{equation}\label{trivbdlelin}
t \mal (x,z) := (t \mal x, \chi(t) z)
\end{equation}

In particular, the character group $\Chara(H)$ canonically embeds
into the group $\Pic_{H}(X)$.
Since we obtain in~\ref{trivbdlelin} indeed any linearization of 
the trivial bundle,
the map $\Chara(H) \to \Pic_{H}(\rq{X})$ 
and the forget map $\Pic_{H}(\rq{X}) \to \Piclin(\rq{X})$
fit together to an exact sequence, compare
also~\cite[Lemma~2.2]{Kr}:
\begin{equation}\label{thesequence}
\xymatrix{
 0  \ar[r] &
{\Chara(H)} \ar[r]^{} &
{\Pic_{H}(\rq{X})} \ar[r] &
{\Piclin(\rq{X})} \ar[r] &
 0 }
\end{equation}

Thus, to obtain $\Piclin(\rq{X}) = 0$,
it suffices to split the map $\Chara(H) \to \Pic_{H}(\rq{X})$
into isomorphisms as follows:
\begin{equation}\label{specialsetting}
\vcenter{%
\xymatrix{%
{\Chara(H)} \ar[rr]^{} \ar[dr]^{{\cong}}_{{\chi^{P} \mapsto P}}& &
{\Pic_{H}(\rq{X})} \\
& {\Pi} \ar[ur]_{q^{*}}^{\cong} &
            }}
\end{equation}

But this is not hard: The fact that $q^{*}$ induces an
isomorphism of $\Pi = \Pic(X)$ and $\Pic_{H}(\rq{X})$ is due 
to~\cite[Prop.~4.2]{Kr}. 
To obtain commutativity, consider $P \in \Pi$.
Choose invertible sections 
$g_i \in \mathcal{A}_{P}(U_{i})$ 
for small open $U_{i}$ covering $X$.
Then the class of $P$ is represented by the bundle 
$P_{\xi}$ arising from the cocycle
\begin{equation}\label{char2pic}
\xi_{ij} 
\; := \; 
\frac{g_{j}}{g_{i}}.
\end{equation}
So the pullback class $q^{*}(P) \in \Pic_{H}(\rq{X})$ is 
represented by the trivially linearized bundle $q^{*}(P_{\xi})$,
which in turn arises from the cocycle
\begin{equation}\label{char2pic2}
q^{*}(\xi_{ij}) 
\; := \; 
q^{*}\left(\frac{g_{j}}{g_{i}}\right)
\; = \; 
\frac{g_{j}}{g_{i}}.
\end{equation}
But on $\rq{X}$, the $g_{i}$ are ordinary invertible functions. 
So we obtain an isomorphism from the representing bundle 
$q^{*}(P_{\xi})$ onto the trivial bundle by locally
multiplying with $g_{i}$. 
Obviously, the induced linearization on the
trivial bundle is the linearization~\ref{trivbdlelin}
for $\chi = \chi^{P}$.

Thus we proved that $(A,\mathfrak{A})$ is in fact a
simple pointed algebra. In other words, the homogeneous coordinate
ring functor restricts to the subcategories in consideration.
It remains to show that up to isomorphism, every simple pointed
algebra is the homogeneous coordinate ring of 
some tame variety $X$.

So, let $(A,\mathfrak{A})$ be a simple pointed algebra, graded by
some finitely generated abelian group $\Pi$. 
According to Proposition~\ref{equivquasiaffequiv}, we may
assume that $(A,\mathfrak{A})$ equals 
$(\mathcal{O}(\rq{X}),\mathfrak{O}(\rq{X}))$ for
some normal quasiaffine variety $\rq{X}$ with a free
action of a diagonalizable group $H = \Spec(\KK[\Pi])$. 

The action of $H$ on $\rq{X}$ admits a geometric quotient 
$q \colon \rq{X} \to X$: First divide by the finite factor 
$H_{{\rm t}}$ of $H$ to obtain a normal quasiaffine variety 
$\rq{X}/H_{{\rm t}}$, and then divide by the induced 
action of the unit component $H_{0}$ of $H$ on 
$\rq{X}/H_{{\rm t}}$,
see for example~\cite[Ex.~4.2]{Do} and~\cite[Cor.~3]{Su}.

The candidate for our tame variety is $X$. Since 
the structure sheaf $\mathcal{O}_{X}$ is the sheaf of invariants 
$q_{*} (\mathcal{O}_{\rq X})^{H}$ and $A = \Of(\rq{X})$ is pointed, we have 
$\mathcal{O}(X) = \KK$. Moreover, as a geometric quotient space of a
normal quasiaffine variety by a free diagonalizable group action, 
$X$ is again normal and divisorial, 
for the latter see~\cite[Lemma~3.3]{Ha1}.
 
To conclude the proof, we have to realize the ($\Pi$-graded) 
direct image $\mathcal{A} := q_{*}(\mathcal{O}_{\rq{X}})$ 
as a Picard graded algebra on $X$. 
First note that we have again the exact 
sequence~\ref{thesequence}.
Since we assumed $\Piclin(\rq{X}) = 0$, the 
character group $\Chara(H)$ maps isomorphically onto
$\Pic_{H}(\rq{X})$.

Moreover, we have a canonical map $\Pi \to \Pic(X)$: For a degree $P \in \Pi$
choose invertible $\chi^{P}$-homogeneous functions $g_{i} \in
\mathcal{O}(q^{-1}(U_{i}))$ with small open $U_{i} \subset X$ 
covering $X$, see Definition~\ref{freegradalgdef}. 
As in~\ref{char2pic}, such functions define a cocycle
$\xi$ and hence we may map $P$ to the class of the bundle $P_{\xi}$.  In
conclusion, we arrive again at a commutative diagram as
in~\ref{specialsetting}.  In particular, $\Pi \to \Pic(X)$ is an isomorphism.

In fact, the construction~\ref{char2pic} allows us to define a
group $\Lambda$ of line bundles on~$X$: 
As in the proof of Lemma~\ref{ontopic}, 
we may adjust the sections $g_{i}$ for a system of 
generators $P$ of $\Pi$,
such that that the corresponding cocycles $\xi$
generate a finitely generated free abelian group. 
Let $\Lambda$ be the resulting group of line bundles, 
and denote the associated $\Lambda$-graded 
$\mathcal{O}_{X}$-algebra by
$\mathcal{R}$.

We construct a graded $\mathcal{O}_{X}$-algebra homomorphsim
$\mathcal{R} \to \mathcal{A}$. The accompanying homomorphism will be
the canonical map $\Lambda \to \Pi$, associating to $L$ its class
under the identification $\Pi \cong \Pic(X)$. Now, the sections of
$\mathcal{R}_{L}$, where $L = P_{\xi}$, 
are given by families $(h_{i})$ satisfying
$$ h_{j} = \xi_{j} h_{i} = \frac{g_{j}}{g_{i}} h_{i} .$$

This enables us to define a map $\mathcal{R}_{L} \to \mathcal{A}_{P}$ by
sending $(h_{i})$ to the section obtained by patching together
the $h_{i}g_{i}$. Note that this indeed yields a graded homomorphism
$\mathcal{R} \to \mathcal{A}$. By construction, this homomorphism is
an isomorphism in every degree. Thus we only have to show that its
kernel is the ideal associated a shifting family for $\mathcal{R}$.

Let $\Lambda_{0} \subset \Lambda$ denote the kernel of the canonical
map $\Lambda \to \Pi$. Then every bundle $E \in \Lambda_{0}$ admits a
global trivialization. In terms of the defining cocycle $g_{i}/g_{j}$
of $E$ this means that there exist invertible local funtions
$\t{g}_{i}$ on $X$ with
$$ 
\frac{g_{j}}{g_{i}} = \frac{\t{g}_{j}}{\t{g}_{i}}.
$$

\goodbreak

The functions $\t{g}_{i}$ can be used to define a shifting
family: Let $L \in \Lambda$ and $E \in \Lambda_{0}$. Then the sections
of $\mathcal{R}_{L}$ are given by families $(h_{i})$ of function that
are compatible with the defining cocycle. Thus we obtain maps
$$ 
\varrho_{E} \colon \mathcal{R}_{L} \to \mathcal{R}_{L+E},
\qquad
(h_{i}) \mapsto \left(\frac{h_{i}}{\t{g}_{i}}\right).
$$

By construction, the $\varrho_{E}$ are homomorphisms, and they
fit together to a shifting family 
$\varrho$ for $\mathcal{R}$. It is straightforward to check that the
ideal $\mathcal{I}$ associated to $\varrho$ is precisely the kernel of
the homomorphism $\mathcal{R} \to \mathcal{A}$. 
\endproof

\section{Very tame varieties}\label{section8}

Finally, we take a closer look to the case of a free Picard group.
The only assumption in this section is that the multiplicative group 
$\KK^{*}$ is of infinite rank over $\ZZ$. 
But even this could be weakened, see 
the concluding Remark~\ref{verytame}.

\begin{defi}
A {\em very tame\/} variety is a normal divisorial variety with
finitely generated free Picard group and only constant functions.
\end{defi}

Examples of very tame varieties are Grassmannians and all smooth
complete toric varieties. On the algebraic side we work with the following 
notion:

\begin{defi}
A {\em very simple\/} algebra is a freely $\Lambda$-graded quasiaffine 
algebra $(A,\mathfrak{A})$ such that
\begin{enumerate}
\item the grading group $\Lambda$ of $(A,\mathfrak{A})$ is free,
\item $A$ is normal, and we have $A_{0} = \KK$ and $A^{*} = \KK^{*}$,
\item the quasiaffine variety associated to $(A,\mathfrak{A})$ 
has trivial Picard group.
\end{enumerate}
\end{defi}

Again, very tame varieties and very simple algebras form subcategories,
and we have an equivalence theorem:

\begin{thm}
The homogeneous coordinate ring functor restricts to an equivalence of
the category of very tame varieties with the category of very simple 
algebras.
\end{thm}

\proof 
Let $X$ be a very tame variety. 
We only have to show is that the quasiaffine $H$-variety 
$\rq{X}$ corresponding to the homogeneous coordinate ring of $X$
has trivial Picard group.
Since $\rq{X}$ is normal and $H$ is a torus, 
every line bundle on $\rq{X}$ is $H$-linearizable,
see~\cite[Remark p.~67]{Kn}.
But from Theorem~\ref{equivthm}, we know that every
$H$-linearizable bundle on $\rq{X}$ is trivial.
\endproof

In the setting of very tame varieties, we can go further with the
dictionary presented in Section~\ref{section6}. The first
remarkable statement is that very tame varieties produce unique
factorization domains:

\begin{prop}\label{freefactorial}
Let $X$ be a very tame variety with homogeneous coordinate ring
$(A, \mathfrak{A})$.
Then $X$ is locally factorial if and only if
$A$ is a unique factorization domain.
\end{prop}

\proof
Let $A = \mathcal{A}(X)$ with some Picard graded
$\mathcal{O}_{X}$-algebra $\mathcal{A}$, and
the geometric quotient $q \colon \rq{X} \to X$ provided by 
Lemma~\ref{geomquot}.
Since $\Pic(X)$ is free we divide by a torus $H$.
Thus $q \colon \rq{X} \to X$ is an $H$-principal bundle
with respect to the Zariski topology.
In particular, $X$ is locally factorial if and only if
$\rq{X}$ is so.
But $\rq{X}$ is locally factorial if and only if
$A$ is a factorial ring, because we have $\Pic(\rq{X}) = 0$.
\endproof

Next we treat products. Let $X$ and $Y$ be very tame varieties
with homogeneos coordinate rings $(A,\mathfrak{A})$ and
$(B,\mathfrak{B})$.
Fix closing subalgebras $(A',I') \in \mathfrak{A}$ and 
$(B',J') \in \mathfrak{B}$,
as in~\ref{freegradalgdef}, and consider the algebra
$$ A \boxtimes B 
:= \bigcap_{f} (A' \otimes_{\KK} B')_{f}
= \bigcap_{f} (A \otimes_{\KK} B)_{f}, $$
where the intersections are taken in the quotient field of $A'
\otimes_{\KK} B'$ and $f$ runs through the elements of the form $g
\otimes h$ with homogeneous $g \in I'$ and $h \in J'$. 

Now $A$ and $B$ are graded, say by $\Lambda$ and $\Gamma$.
These gradings give rise to a  
$(\Lambda \times \Gamma)$-grading of $A \boxtimes B$. 
Moreover,
$$(A' \otimes_{\KK} B', \sqrt{I' \otimes_{\KK} J'})$$
is a closing subalgebra of $A \boxtimes B$. Let $\mathfrak{A}
\boxtimes \mathfrak{B}$ denote the equivalence class of this closing
subalgebra. Then we obtain:

\begin{prop}\label{products}
Let $X$ and $Y$ be locally factorial very tame varieties.
Then $X \times Y$ is locally factorial and very tame with 
homogeneous coordinate ring $(A \boxtimes B, \mathfrak{A}
\boxtimes \mathfrak{B})$. 
Moreover, if $A$ and $B$ are of finite
type over $\KK$, then $A \boxtimes B$ equals $A \otimes_{\KK} B$.
\end{prop}

\proof 
First note that for any two quasiaffine varieties 
$\rq{X}$ and $\rq{Y}$ with free diagonalizable group 
actions, their product $\rq{X} \times \rq{Y}$ is again 
such a variety.
Moreover, if $\rq{X}$ and $\rq{Y}$ have only constant 
invertible functions, then so does $\rq{X} \times \rq{Y}$.
If $\rq{X}$ and $\rq{Y}$ are additionally 
locally factorial with trivial Picard groups,
then the same holds for $\rq{X} \times \rq{Y}$, 
use e.g.~\cite[Prop.~1.1]{FI}.

Now, let $\rq{X} := \Spec(\mathcal{A})$ and 
$\rq{Y} := \Spec(\mathcal{B})$. 
By Proposition~\ref{freefactorial} both are locally factorial.
By construction
$(A \boxtimes B, \mathfrak{A} \boxtimes \mathfrak{B})$ 
is the freely graded quasiaffine algebra corresponding 
to the product $\rq{X} \times \rq{Y}$.
Thus the above observations and
Proposition~\ref{equivquasiaffequiv} tell us that
it is a coproduct in the category of simple pointed algebras. 
Hence the assertion follows from Theorem~\ref{equivthm}. 
The second statement is an
easy consequence of Remark~\ref{collaps}~(i).
\endproof 

\begin{coro}
Let $X$ and $Y$ be locally factorial very tame varieties. 
Then $\Pic(X \times Y)$ is isomorphic to 
$\Pic(X) \times \Pic(Y)$. \endproof
\end{coro}

We give an explicit example
emphasizing the role of Proposition~\ref{freefactorial}. 
We assume that the ground field $\KK$ is not of characteristic two. 
Consider the prevariety $X$ obtained by gluing two copies of 
the projective line 
$\PP_{1}$ along the common open subset $\KK^{*} \setminus \{1\}$. 
We think of $X$ as the projective line with three doubled points,
namely
$$ 0, \; 0', \quad 1,  \; 1', \quad \infty, \; \infty'.$$

Note that $X$ is smooth and divisorial. Moreover, $\Pic(X)$ is
isomorphic to $\ZZ^{4}$. Thus we obtain in particular that $X$ is
very tame. Let $(\mathcal{A}(X),\mathfrak{A}(X))$ denote the homogeneous
coordinate ring of $X$. We show:

\begin{prop}\label{example}
$\mathcal{A}(X) 
\cong 
\KK[T_{1}, \ldots, T_{6}]/\bangle{T_{1}^{2} + \ldots + T_{6}^{2}}$.
\end{prop}

Before giving the proof, let us remark that the ring $\mathcal{A}(X)$
is a classical example of a singular factorial affine algebra.
In view of our results, factoriality is a consequence of
Proposition~\ref{freefactorial}.

\proof[Proof of Proposition~\ref{example}]
First observe that we may realize $\Pic(X)$ as well as a subgroup 
$\Lambda$ of the group of Cartier divisors of $X$. For example
$\Pic(X)$ is isomorphic to the group $\Lambda$ generated by 
$$
D_{0} := \{0\},
\qquad
D_{1} := \{1\},
\qquad
D_{1'} := \{1'\},
\qquad
D_{\infty} := \{\infty\}.$$

For any Cartier divisor $D$ on $X$, let $\mathcal{A}_{D}$ denote its
sheaf of sections. Then the homogeneous coordinate ring
$\mathcal{A}(X)$ is the direct sum of the $\mathcal{A}_{D}(X)$, 
where $D \in \Lambda$. Consider the following homogeneous elements 
of $\mathcal{A}(X)$: 
$$ 
\begin{array}{ll}
f_{1} := 1 \in \mathcal{A}_{D_{0}} (X), 
&
f_{2} := 1 \in \mathcal{A}_{D_{1}} (X), \cr
f_{3} := 1 \in \mathcal{A}_{D_{1'}}(X),
&
f_{4} := 1 \in \mathcal{A}_{D_{\infty}}(X). \cr
f_{5} := \bigl( \frac{1}{z-1} \bigr) 
   \in \mathcal{A}_{D_{1} + D_{1'} - D_{\infty}}(X),
&
f_{6} := \bigl( \frac{z}{z-1} \bigr)
 \in \mathcal{A}_{D_{1} + D_{1'} - D_{0}}(X).
\end{array}
$$

Let $\varphi$ be the algebra homomorphism $\KK[T_{1}, \dots, T_{6}] \to
\mathcal{A}(X)$ sending $T_{i}$ to $f_{i}$. It is elementary to check
that $\varphi$ is surjective. Since we assumed $\KK$ not to be of
characteristic two, it suffices to show that the kernel of $\varphi$ 
is the ideal generated by
$$
Q := T_{2}T_{3} + T_{5} T_{4} - T_{6} T_{1}.
$$

An explicit calculation shows that the $f_{i}$ 
fulfil the claimed relation, that means that $Q$
lies in the kernel of $\varphi$. 
Conversely, consider an arbitrary element $R$ of the kernel of $\varphi$.
Then there are $r_{j} \in \KK[T_{1}, \dots, T_{5}]$ such that
$R$ is of the form
$$
R = \sum_{j=0}^{s} r_{j} T_{6}^{j}.
$$

We proceed by induction on $s$. 
For $s = 0$ the fact that $f_{1}, \dots, f_{5}$ are algebraically
independent implies $R = 0$.
For $s > 0$ note first that that $\varphi(r_{j})$
is nonnegative in $D_{0}$ in the sense that
its component in a degree containing a multiple 
$nD_{0}$ is  trivial for negative $n$.
 
Since $f_{6}$ is negative in $D_{0}$, and $f_{1}$ is 
the only generator of $\mathcal{A}(X)$
which is strictly positive in degree $D_{0}$,
we can write $r_{j} = \t{r}_{j} T_{1}^{j}$. 
Hence we obtain a representation
$$ R = \sum_{j=0}^{s} \t{r}_{j} T_{1}^{j}T_{6}^{j}.$$
The element 
$\t{r}_{s} ((T_{1}T_{6})^{s} - (T_{2}T_{3} + T_{4}T_{5})^{s})$
is a multiple of~$Q$. 
In particular, it belongs to the kernel of $\varphi$.
Subtracting it from $R$, we obtain 
$$ R'= \sum_{j=0}^{s-1} r_{j}' T_{1}^{j} T_{6}^{j},$$ 
with $r_{j}' = \t{r}_{j}$ for $j>0$ and $r_{0}' = \t{r}_{0} +
\t{r}_{s}(T_{2}T_{3} + T_{4}T_{5})^{j}$. 
Applying the induction hypothesis to $R'$ 
yields that $R$ is a multiple of~$Q$.
\endproof

Finally, let us note that all our statements on very tame varieties
hold under more general assumptions. 
This is due to the fact 
that free Picard groups always can be realized by (free)
groups of line bundles. 
Hence in this case we don't need shifting families to define
the homogeneous coordinate ring. 
This means:

\begin{rem}\label{verytame}
For very tame varieties $X$, the results of this article
hold over any algebraically closed ground field $\KK$, 
and one might weaken the assumption $\mathcal{O}(X) = \KK$ 
to $\mathcal{O}^{*}(X) = \KK^{*}$. 
\end{rem}

\end{document}